\RequirePackage{fix-cm}
\documentclass[smallextended]{svjour3}       
\smartqed  
\usepackage{graphicx}
\usepackage{url}
\usepackage{amssymb,amsmath,a4wide,mdwlist}
\usepackage{todonotes}
\usepackage{caption}
\usepackage{subcaption}
\usepackage{blkarray}
\usepackage{hyperref}
\usepackage[linesnumbered,lined,boxed,commentsnumbered]{algorithm2e}
\usepackage{algorithmic}
\newenvironment{proof1}{\noindent {\bf Proof. }}{\hfill $\Box$ \newline\par}

\begin{document}

\title{Computing mixed strategies equilibria in presence of switching costs by the solution of nonconvex QP problems 
}


\author{G. Liuzzi         \and
        M. Locatelli \and V. Piccialli \and S. Rass 
}


\institute{G. Liuzzi \at
             Istituto di Analisi dei Sistemi ed Informatica "Antonio Ruberti" (IASI)    \\
Consiglio Nazionale delle Ricerche (CNR)        \\
Via dei Taurini 19, 00185 Rome - Italy \\
              Tel.: +0039 06 4993 7129\\
              \email{giampaolo.liuzzi@cnr.iasi.it}           
           \and
           M. Locatelli \at
              Universit\`{a}  degli Studi di Parma    \\
                Parco Area delle Scienze, 181/A - I 43124 Parma     \\\email{marco.locatelli@unipr.it}
                \and
                V. Piccialli \at
                DICII - University of Rome Tor Vergata \\
                via del Politecnico 1 \\
                00133 Roma\\
                \email{veronica.piccialli@uniroma2.it}
                \and
                S. Rass \at
                Universit\"at Klagenfurt, Institute of Applied Informatics, System Security Group,\\
Klagenfurt, Austria, \\
\email{stefan.rass@aau.at}
}

\date{Received: date / Accepted: date}

\maketitle

\begin{abstract}
In this paper we address game theory problems arising in the context of network security. In traditional game theory problems, given a defender and an attacker, one searches for mixed strategies  which minimize a linear payoff functional. 
In the problems addressed in this paper an additional quadratic term is added to the minimization problem. Such term represents {\em switching costs}, i.e., the costs for the defender of switching from a given strategy to another one at successive rounds of a Nash game. The resulting problems are nonconvex QP ones with linear constraints and
turn out to be very challenging. We will show that the most recent approaches for the minimization of nonconvex QP functions over polytopes, including commercial solvers such as {\tt CPLEX} and {\tt GUROBI}, are unable to solve to optimality even test instances with $n=50$ variables. For this reason, we propose to extend with them the current benchmark set of test instances for QP problems.
We also present a spatial branch-and-bound approach for the solution of these problems,
where a predominant role 
is played by an optimality-based domain reduction, with multiple solutions of LP problems at each node of the branch-and-bound tree. Of course, domain reductions are standard tools in spatial branch-and-bound approaches. However, our contribution lies in the observation that, from the computational point of view, a rather aggressive application of these tools appears to be the best way to tackle the proposed instances. Indeed, according to our experiments, while they make the computational cost per node high, this is largely compensated by the rather slow growth of the number of nodes in the branch-and-bound tree, so that the proposed approach strongly outperforms the existing solvers for
QP problems.
\end{abstract}

\section{Introduction}
\label{intro}

Consider a finite two-person zero-sum game $\Gamma$, composed from a player
set $N=\{1,2\}$, each member thereof having a finite strategy space $S_1,S_2$
associated with it, and a utility function $u_i:S_1\times S_2\to\mathbb{R}$
for all $i\in N$. We assume a zero-sum Nash game, making $u_2:=-u_1$
hereafter, and letting the players choose their actions simultaneously and stochastically independent of one another (contrary to a Stackelberg game, where one player would follow the other, which we do not consider here). The game is then the triple $\Gamma=(N,\mathcal
S=\{S_1,S_2\},H=\{u_1,-u_1\})$, and is most compactly represented by giving
only the payoff function $u_1$ in matrix form (since the strategy spaces are
finite) as
\[
    {\bf A}\in\mathbb{R}^{|S_1|\times|S_2|}=\left(u_1(x,z)\right)_{(x,z)\in S_1\times S_2}.
\]

An equilibrium in $\Gamma$ is a simultaneous optimum for both players w.r.t.
$u_1$. Assuming a maximizing first player, an equilibrium is a pair
$(x^*,z^*)$ satisfying the saddle-point condition
\[
    u_1(x,z^*)\leq u_1(x^*,z^*)\leq u_1(x^*,z)\quad\forall (x,z)\in S_1\times S_2.
\]
It is well known that many practical games do not have such an equilibrium
point; as one of the simplest instances, consider the classical
rock-scissors-paper game, represented by the payoff matrix
\[
\begin{blockarray}{rccc}
& \text{rock} & \text{ scissors } & \text{ paper }  \\
\begin{block}{r(ccc)}
  \text{ rock } & 0 & 1 & -1  \\
  \text{ scissors } & -1 & 0 & 1  \\
  \text{ paper } & 1 & -1 & 0  \\
\end{block}
\end{blockarray}.
\]
This game has no equilibria in pure strategies: any fixed choice of rock, scissors or
paper would imply a constant loss for the first player (and likewise for the
second player). This means that player 1 is forced to \emph{randomize} its
actions in every round of the game, and this concept leads to the idea of
mixed extensions of a game, which basically changes the above optimization
problem into one over the convex hulls $\Delta(S_1), \Delta(S_2)$ of the
action spaces, rather than the finite sets $S_1,S_2$. An element of
$\Delta(S_i)$ is then a probability distribution over the elements of the
support $S_i$, and prescribes to pick a move at random whenever the game is
played.
\newline\newline\noindent
The game rewards its players after each round, and upon every new round, both
players are free to choose another element from their action space at random.
Implicitly, this choice is without costs, but what if not? Many real life
instances of games \emph{do incur} a cost for changing one's action from
$a_1\in S_1$ in the first to some distinct $a_2\in S_1$ in the next round. Matrix games cannot express such costs in their payoff functions, and more complex game models such as sequential or stochastic games come with much more complicated models and equilibrium concepts. The goal of this work is to retain the simplicity of matrix games but endow them with the ability to include switching costs with the minimal natural (modeling) effort. 

The
area of system security \cite{Alpcan&Basar2010,Tambe2012} offers rich
examples of such instances, such as (among many):
\begin{itemize}
  \item Changing passwords \cite{rass_password_2018}: if the currently
      chosen password is $p_1$ and we are obliged to pick a fresh password
      (say, different from the last couple of passwords that we had in the
      past), the use of the new password $p_2\neq p_1$ induces quite some
      efforts, as we have to memorize the password, while choosing it as
      hard as possible to guess. The ``cost'' tied to the change is thus not monetary, but
      the cognitive efforts to create and memorize a new password. This effort can 
      make people reluctant to change their passwords (or write them down, or use a very 
      similar password for the new one).
  \item Changing computer/server configurations: this usually means taking a computer (e.g., a server) offline for a limited time, thus cutting down productivity perhaps, and hence causing costs. If security is drawn from randomly changing configurations (and passwords, resp. password changing rules are only one special case here), then this change incurs costs by temporal outages of IT infrastructure for the duration of the configuration change, and the efforts (person-hours) spent on applying this change. This is why server updates or patches are usually done over nights or weekends, when the loads are naturally low. If the optimization would, however, prescribe a rather frequent change of configurations at random intervals, this can quickly become a practical inhibitor, unless the switching costs are accounted for by optimization.
  \item Patrolling and surveillance \cite{Alpern.2011,rass_physical_2017}:
      consider a security guard on duty to repeatedly check a few
      locations, say A, B, \ldots, E, which are connected at distances as
      depicted in Fig. \ref{fig:graph-example}.
\begin{figure}[h!]
  \centering
  \includegraphics[scale=1]{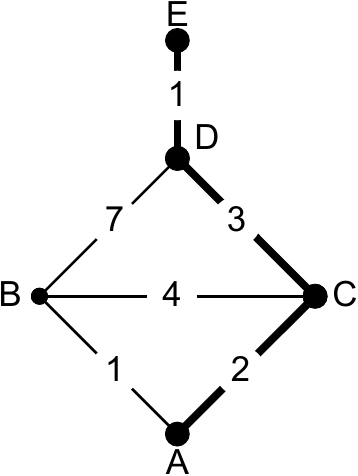}
  \caption{Example of spot checking game on a graph}\label{fig:graph-example}
\end{figure}
This is a chasing-evading game with the guard acting as player 1
against an intruder being player 2, and with the payoff function $u$ being an indicator of whether the guard caught the intruder at location $i\in\{$A,\ldots,E$\}$, or whether the two missed each
other. This is yet another instance of a game with all equilibria in
mixed strategies, but with the unpleasant side-effect for the guard that gets the prescription to randomly spot check distant locations to
``play'' the equilibrium ${\bf x^*}$, the guard would have to move
perhaps long distances between the locations. For example, if it is
at A in round 1 and the next sample from the random distribution
${\bf x^*}\in\Delta(\{$A,\ldots,E$\})$ tells to check point E next,
the shortest path would be of length $1+3+2=6$ over C. Starting from
A, however, it would be shorter and hence more convenient for the
guard to check location B first along the way, but this would mean
deviating from the equilibrium! A normal game theoretic equilibrium
calculation does not consider this kind of investment to change the
current strategy. This may not even count as bounded rationality, but
simply as acting ``economic'' from the guard's perspective. But
acting economically here is then not governed by a utility maximizing
principle, but rather by a cost minimization effort.
\item Generalizing the patrolling game example, the issue applies to all
sorts of moving target defense: for example, changing the
configuration of a computer system so as to make it difficult for an
attacker to break in, often comes with high efforts and even risks
for the defending player 1 (the system administrator), since it
typically means taking off machines from the network, reconfiguring
them to close certain vulnerabilities, and then putting them back to
work hoping that everything restarts and runs smoothly again. A
normal game theoretic model accounts only for the benefits of that
action, but not for the cost of \emph{taking} the action.
\end{itemize}
Including the cost to switch from one action to the next is more complicated
than just assigning a cost function $c:S_1\to\mathbb{R}$ and subtracting this
from the utilities to redefine them as $u_1'(i,j)=u_1(i,j)-c(i)$, since the
cost to play $a_i$ will generally depend on the previous action $a_j$ played
in the previous round.
\newline\newline\noindent
We can model this sort of payoff by another function $s:S_1\times
S_1\to\mathbb{R}$ that we call the \emph{switching cost}. The value of
$s(i,j)$ is then precisely the cost incurred to change the current action
$i\in S_1$ into the action $j\in S_1$ in the next round of the game.
Intuitively, this adds another payoff dimension to the game, where a
player, w.l.o.g. being player 1 in the following, plays ``against itself'',
since the losses are implied by its own behavior. While the expected payoffs
in a matrix game ${\bf A}$ under mixed strategies ${\bf x}\in\Delta(S_1),{\bf
z}\in\Delta(S_2)$ are expressible by the bilinear functional ${\bf x}^T{\bf
A}{\bf z}$, the same logic leads to the hypothesis that the switching cost
should on average be given by the quadratic functional ${\bf x}^T{\bf S}{\bf
x}$, where the switching cost matrix is given, like the payoff matrix above,
as
\[
    {\bf S}\in\mathbb{R}^{|S_1|\times|S_1|}=\left(s(x,w)\right)_{(x,w)\in S_1\times S_1}.
\]
This intuition is indeed right \cite{Rass.2017}, but for a rigorous problem
statement, we will briefly recap the derivation given independently later by
\cite{Wachter.2018} to formally state the problem.
\subsection{Paper Outline}
The paper is structured as follows. In Section \ref{sec:1} we give a formal description of the problem as a nonconvex QP one with linear constraints, and we report a complexity result, proved in Appendix
\ref{sec:complexity}.
In Section \ref{sec:bbsec} we present a (spatial) branch-and-bound
approach for the problem, putting a particular emphasis on the bound-tightening procedure,
which turns out to
be the most effective tool to attack it. In Section \ref{sec:numrel} we present some computational
experiments. We first describe the set of test instances. Next, we discuss the performance of existing
solvers over these instances. Finally, we present and comment the computational results attained
by the proposed approach.
In Section \ref{sec:concl} we draw some conclusions and discuss possible future developments.
\subsection{Statement of contribution}
The main contributions of this work are:
\begin{itemize}
\item addressing an application of game theory arising in the context of network security, where switching costs come into play, and showing that the resulting problem can be reformulated as a challenging nonconvex QP problem with linear constraints;
\item introducing a large set of test instances, which turn out to be very challenging for existing QP solvers and, for this reason, could be employed to extend the current benchmark set of QP problems (see \cite{furini2019qplib});
\item proposing a branch-and-bound approach for the solution of the addressed QP problems, based on standard tools, but with the empirical observation that a very aggressive use of bound-tightening techniques, with a high computational cost per node of the branch-and-bound tree, is the key for an efficient solution of these problems.
\end{itemize}
\section{Formal description of the problem}\label{sec:1} 
Let the game come as a matrix ${\bf A}\in\mathbb{R}^{n\times
m}$, where $n$ and $m$ are the number of strategies for player 1 and 2, respectively, with equilibrium $({\bf x}^*,{\bf z}^*)$, and let it be repeated over the
time $t\in\mathbb{N}$. At each time $t$, let $X_t\sim {\bf x}^*$ be the
random action sampled from the equilibrium distribution over the action space
(with ${\bf x}^*$ being the optimal distribution). In a security setting and zero-sum game, neither player has an interest of being predictable by its opponent, so we assume stochastic independence of the action choices by both players between any two repetitions of the game. Then, we have $\Pr(X_{t-1}=i,X_t=j)=\Pr(X_{t-1}=i)\cdot\Pr(X_{t}=j)$, so that
any future system state remains equally predictable whether or not the
current state of the system is known. Hence, the switching cost can be
written as
\begin{align*}
s(X_{t-1}, X_t)&=\sum_{i=1}^{n}\sum_{j=1}^{n}s_{ij}\cdot\Pr(X_{t-1}=i,X_t=j)\\
&=\sum_{i=1}^{n}\sum_{j=1}^{n}s_{ij}\cdot\Pr(X_{t-1}=i)\cdot\Pr(X_{t}=j) = {\bf
x}^T{\bf S}{\bf x}.
\end{align*}
With this, player 1's payoff functional becomes vector-valued now as
\begin{equation}\label{eqn:multi-objective-payoff}
    {\bf u}_1:\Delta(S_1)\times \Delta(S_2)\to\mathbb{R}^2,\quad
    ({\bf x},{\bf z})\mapsto \left(
                               \begin{array}{l}
                                 u_1({\bf x},{\bf z})={\bf x}^T{\bf A}{\bf z} \\
                                 s({\bf x},{\bf z})={\bf x}^T{\bf S}{\bf x} \\
                               \end{array}
                             \right),
\end{equation}
and the game is multi-objective for the first player. As we are interested
mostly in the best behavior for player 1 and the analysis would be symmetric
from player 2's perspective, we shall not explore the view of the second
player hereafter.
\begin{remark}The game could
be equally well multi-objective for the second player too, and in fact a
practical instance of such a situation may also come from security: it could
be in an adversary's interest to ``keep the defender busy'', thus causing
much friction by making the defender move fast from one place to the other.
This is yet just another instance of a denial-of-service attack, to which
such a game model would apply.
\end{remark}
$\ $\newline\newline\noindent
For the sake of computing a multi-objective equilibrium, more precisely a
Pareto-Nash equilibrium, the algorithm in \cite{rass_numerical_2014} based on
the method laid out in \cite{lozovanu_multiobjective_2005} proceeds by
scalarizing \eqref{eqn:multi-objective-payoff} by choice of some
$\alpha\in(0,1)$, to arrive at the real-valued goal function
\[
    \alpha\cdot{\bf x}^T{\bf A}{\bf z}+(1-\alpha)\cdot{\bf x}^T{\bf S}{\bf x},
\]
for the first player to optimize. Now, the usual way from here to an
optimization problem for player 1 involving a rational opponent applies as
for standard matrix games \cite{Rass.2017}: let ${\bf e}_i\in\mathbb{R}^m$ be
$i$-th unit vector, then $\arg\max_{{\bf z}\in \Delta(S_2)}({\bf x}^T{\bf
A}{\bf z})=\arg\max_{i}({\bf x}^T{\bf A}{\bf e}_i)$. After introducing the additional variable $v$, the resulting problem becomes
\begin{equation}\label{eqn:nlp}
      \begin{array}{llll}
       \min & &(1-\alpha)\cdot{\bf x}^T{\bf S}{\bf x}+ \alpha v \\[6pt]
        &\text{s.t.} & v  \geq {\bf x}^T{\bf A}{\bf e}_i & i=1,\ldots,m \\ [6pt]
         &                        & \sum_{j=1}^n x_j  =  1 & \\[6pt]
        & & x_j \geq 0 & j=1,\ldots,n,
      \end{array}
\end{equation}
which is almost the familiar optimization problem to be solved for a Nash
equilibrium in a finite matrix game. It differs from the well known linear
program only in the quadratic term, and, in fact, the equilibrium problem for
matrix games is recovered by substituting $\alpha=1$ in \eqref{eqn:nlp}. The
problem would again become trivial for $\alpha=0$, since in that case, only
the switching cost matters and hence every degenerate distribution
corresponding to a strategy $i\in S_1$ that never changes is directly an
equilibrium. Excluding the standard equilibria obtained at $\alpha=1$ and the
meaningless results expected for $\alpha=0$, problem \eqref{eqn:nlp} is
interesting only for values of $\alpha$ strictly between 0 and 1. 
Note that the matrix ${\bf S}$ in the quadratic term will (in most cases) have a zero diagonal, nonnegative off-diagonal entries, be indefinite and not symmetric in general (patrolling game example given above already
exhibits a variety of counterexamples leading to nonsymmetric distance matrices ${\bf S}$ if the graph is directed). Of course, symmetry of ${\bf S}$ can be easily recovered, so in what follows we will assume that ${\bf S}$ is symmetric. 
As already commented, the two extreme values $\alpha=0$ and $\alpha=1$ give rise to simple problems. Indeed, for $\alpha=1$ the problem is an LP one, while for $\alpha=0$ is a Standard QP (StQP) problem, which is in general NP-hard (e.g., in view of the reformulation of the max clique problem as an StQP problem, see \cite{Motzkin65}), but is trivial in the  case of zero diagonal and nonnegative off-diagonal entries (each vertex of the unit simplex is a globally optimal solution).
For what concerns the intermediate values $\alpha\in (0,1)$ we can prove the following result, stating the complexity of problem (\ref{eqn:nlp}) .
\begin{theorem}
\label{theo:compl}
Problem (\ref{eqn:nlp}) is NP-hard.
\end{theorem}
\begin{proof1}
See Appendix \ref{sec:complexity}.
\end{proof1}

\begin{remark}
Observe the interesting effect that the two extreme instances at $\alpha=0$ and $\alpha=1$ are solvable in polynomial time, while any intermediate instance with $0<\alpha<1$ is NP-complete. The jump in the complexity thus cannot be attributed to either term alone, but only to their coincidental presence.
\end{remark}

\begin{remark}
The dependence of next actions on past ones extends to other scenarios too:
for example, if the game is about coordination in wireless settings (e.g.,
collaborative drones), the players, e.g., drones, share a common
communication channel. Every exchange of information occupies that channel
for a limited period of time, thus constraining what the other players can do
at the moment. Such effects can be described by stochastic games, but
depending on how far the effect reaches in the future, backward inductive
solution methods may become computationally infeasible \cite{Hansen2012};
likewise, extending the strategy space to plan ahead a fixed number of $k$
steps (to account for one strategy determining the next $k$ repetitions of
the game) may exponentially enlarge the strategy space (by a factor of
$2^{O(k)}$, making the game infeasible to analyze if $k$ is large). Games
with switching cost offer a neat bypass to that trouble: if an action is such
that it occupies lots of resources for a player, thus preventing it from
taking further moves in the next round of the game, we can express this as a
switching cost. Assume, for instance, that an action in a game $\Gamma$ is
such that the player is blocked for the next $k$ rounds, then the switching
cost is $k$-times the expected utility $\overline{u}$ (with the expectation
taken over the equilibrium distribution played by the participants) that
these next $k$ rounds would give. Virtually, the situation is thus like if
the player would have paid the total average gain over the next rounds where
it is forced to remain idle (thus gaining nothing):
\begin{equation}\label{eqn:virtual-cost}
  \overline{u} \underbrace{- k\cdot \overline{u}}_{\text{switching cost}} +
  \underbrace{\overline{u} + \cdots + \overline{u}}_{\text{\shortstack{virtual payoffs\\over $k$ rounds}}} = \overline{u} + \underbrace{0 + 0 + \ldots + 0}_{\text{\shortstack{practical payoffs\\by being idle\\for $k$ rounds}}}
\end{equation}
$\ $\newline\newline\noindent
Expression \eqref{eqn:virtual-cost} will in practice be only an approximate
identity, since we assumed that the game, viewed as a stochastic process, has
already converged to stationarity (so that the equilibrium outcome
$\overline{u}$ is actually rewarded). The speed of convergence, indeed, can
itself be of interest to be controlled in security applications using moving
target defenses \cite{Wachter.2018}. The crucial point of modeling a longer
lasting effect of the current action like described above, however, lies in
the avoidance of complexity: expression \eqref{eqn:virtual-cost} has no
issues with large $k$, while more direct methods of modeling a game over $k$
rounds, or including a dependency on the last $k$ moves, is relatively more
involved (indeed, normal stochastic games consider a first-order Markov
chain, where the next state of the game depends on the last state; the
setting just described would correspond to an order $k$ chain, whose
conversion into a first order chain is also possible, but complicates matters
significantly).
\end{remark}

\section{A branch and bound approach}
\label{sec:bbsec}
After incorporating parameter $\alpha$ 
into the definitions of matrix ${\bf S}$ and vectors ${\bf A}_j$, $j=1,\ldots,m$, and after introducing the vector of variables ${\bf y}$, problem (\ref{eqn:nlp}) can be rewritten as the following problem 
with bilinear objective function and linear constraints:
\begin{equation}
\label{eq:bilref}
\begin{array}{lll}
\min & F({\bf x}, {\bf y}, v):= \frac{1}{2}\sum_{i=1}^n x_i y_i  + v & \\ [8pt]
& v \geq {\bf A}_j^T {\bf x}& j=1,\ldots,m \\ [8pt]
& y_i={\bf S}_i {\bf x} & i=1,\ldots,n \\ [8pt]
& {\bf x}\in \Delta_n,  &
\end{array}
\end{equation}
where ${\bf S}_i$ denotes the $i$-th row of matrix ${\bf S}$ and $\Delta_n$ denotes the $n$-dimensional unit simplex.
In what follows we will denote by $P$ the feasible region of this problem, and by $P_{{\bf x},{\bf y}}$ its projection over the space of ${\bf x}$ and ${\bf y}$ variables.
\newline\newline\noindent
Each node of the branch-and-bound tree is associated to a box $B=[{\bf \ell}_{{\bf x}}, {\bf u}_{{\bf x}}]\times [{\bf \ell}_{{\bf y}},{\bf u}_{{\bf y}}]$, where ${\bf \ell}_{{\bf x}}, {\bf u}_{{\bf x}}$ and ${\bf \ell}_{{\bf y}}, {\bf u}_{{\bf y}}$
denote lower and upper bound vectors for variables ${\bf x}$ and ${\bf y}$, respectively.
An initial box $B_0$, containing $P_{{\bf x},{\bf y}}$ is easily computed. It is enough to set ${\bf \ell}_{{\bf x}}={\bf 0}$, ${\bf u}_{{\bf x}}={\bf e}$ (the vector whose entries are all equal to one), and
$$
\ell_{y_i}=\min_{k=1,\ldots,n} S_{ik},\ \ \ \ell_{u_i}=\max_{k=1,\ldots,n} S_{ik}.
$$
Note that, although not strictly necessary, we can also bound variable $v$ to belong to an interval. Indeed, we can impose
$v\geq 0$ (due to nonnegativity of the entries of vectors ${\bf A}_j$, $j=1,\ldots,m$), and 
$$
v\leq \max_{j=1,\ldots,m,\ k=1,\ldots,n} A_{jk},
$$
which certainly holds at optimal solutions of problem (\ref{eq:bilref}).
In what follows we describe in detail each component of the branch-and-bound approach, whose pseudo-code is then sketched in Algorithm \ref{Alg:BandB}.
\subsection{Lower bounds} 
Given box $B=[{\bf \ell}_{{\bf x}}, {\bf u}_{{\bf x}}]\times [{\bf \ell}_{{\bf y}},{\bf u}_{{\bf y}}]$, 
then the well known McCormick underestimating function (see \cite{McCormick76})
$$
\max\left\{\ell_{x_i} y_i+\ell_{y_i} x_i -\ell_{x_i} \ell_{y_i},  u_{x_i} y_i+u_{y_i} x_i -u_{x_i} u_{y_i} \right\},
$$
can be employed to limit from below the bilinear term $x_i y_i$ over
the rectangle $[\ell_{x_i}, u_{x_i}]\times [\ell_{y_i}, u_{y_i}]$.
In fact, it turns out that McCormick underestimating function is the convex envelope of the bilinear term over the given rectangle.
Then, after introducing the additional variables $f_i$, we have that
the optimal value of the following LP problem is a lower bound for problem (\ref{eq:bilref}) over the box $B$:
\begin{subequations}\label{boundLP}
\begin{eqnarray}\label{obj}
&L(B)=\min & \frac{1}{2}\sum_{i=1}^n f_i+v\\[8pt]\label{simp}
&&{\bf x}\in \Delta_n \\[8pt]\label{nneg}
&&v\geq {\bf A}_j^T {\bf x}\quad j=1,\ldots,m\\[8pt]\label{branch1}
&& y_i={\bf S}_i {\bf x} \quad i = 1,\ldots,n \\[8pt]\label{ydef}
&& ({\bf x}, {\bf y})\in B \\[8pt]\label{branch2}
&& f_i \geq \ell_{y_i}x_i + \ell_{x_i} y_i -\ell_{x_i} \ell_{y_i} \quad i=1,\ldots,n\\[8pt]\label{MC2}
&& f_i\geq u_{x_i} y_i+u_{y_i}x_i-u_{y_i}u_{x_i}\quad i=1,\ldots,n.
\end{eqnarray}\end{subequations}
The optimal solution of the LP problem will be denoted by $({\bf x}^\star(B), {\bf y}^\star(B),{\bf f}^\star(B), v^\star(B))$.
\subsection{Upper bound}
The global upper bound (GUB in what follows) can be initialized with $+\infty$ or, alternatively, if a local search procedure is available, one may run a few local searches from randomly generated starting points, and take the lowest local minimum value as initial GUB value, although, according to our experiments, there is not a significant variation in the computing times if such local searches are performed. During the execution of the branch-and-bound algorithm, each time we compute the lower bound
(\ref{boundLP}) over a box $B$, its optimal solution is a feasible solution for problem (\ref{eq:bilref}) and, thus, we might update the upper bound as follows:
$$
GUB=\min\{GUB, F({\bf x}^\star(B), {\bf y}^\star(B), v^\star(B))\}.
$$
\subsection{Branching}
\label{sec:branching}
The branching strategy we employed is a rather standard one. Given node 
$B$, we first compute the quantities:
\begin{equation}\label{gap_i}
  g_i=  x_i^\star (B) y_i^\star(B)-f_i^\star(B),
\end{equation}
measuring the error of McCormick underestimator for each bilinear term $x_i y_i$ at the optimal solution of the relaxed problem (\ref{boundLP}).
Then, we select $r\in\arg\max_{i=1,\ldots,n} g_i$, i.e., the index corresponding to the bilinear term where we have the largest error at the optimal solution of the relaxation. Next, we might define the following branching operations for box $B$:
\begin{description}
\item[{\bf Branching on $x$ and $y$:}] Define four children nodes by adding constraints $\{x_r\leq x_r^\star(B),\ y_r\leq y_r^\star(B)\}$, $\{x_r\leq x_r^\star(B),\ y_r\geq y_r^\star(B)\}$, $\{x_r\geq x_r^\star(B),\ y_r\leq y_r^\star(B)\}$,  $\{x_r\geq x_r^\star(B),\ y_r\geq y_r^\star(B)\}$, respectively (quaternary branching);
\item[{\bf Branching on $x$:}] Define two children nodes by adding constraints $x_r\leq x_r^\star(B)$ and $x_r\geq x_r^\star(B)$, respectively (binary branching);
\item[{\bf Branching on $y$:}] Define two children nodes by adding constraints $y_r\leq y_r^\star(B)$ and $y_r\geq y_r^\star(B)$, respectively (binary branching).
\end{description}
Note that all choices above, with the new McCormick relaxation given by the new limits on the variables, reduce to zero the error for bilinear term $x_r y_r$ at the optimal solution of problem (\ref{boundLP}).
It is worthwhile to remark that the computed lower bound tends to become exact even
when branching is always performed with respect to variables of the same type (say, always variables $x_i$, $i=1,\ldots,n$).
Indeed, it is enough to have that  $\|{\bf u}_{{\bf x}}-{\bf \ell}_{{\bf x}}\| \rightarrow 0$ or, alternatively, 
that $\|{\bf u}_{{\bf y}}-{\bf \ell}_{{\bf y}}\| \rightarrow 0$ in order to let the underestimating function values converge to the original objective function values. This is a consequence of the fact that the McCormick underestimation function tends to the value of the corresponding bilinear term even when only one of the two intervals on which it is defined shrinks to a single point.
In the computational experiments we tried all three possibilities discussed above and it turns out that
the best  choice is the binary branching obtained by always branching on $y$ variables.
\subsection{Bound-tightening technique}\label{sec:boundtight}
A reduction of the boxes merely based on the above branching strategy would lead to a quite inefficient algorithm. It turns out that performance can be strongly enhanced by an Optimality-Based Bound-Tightening (OBBT in what follows) procedure (see, e.g., \cite{DomRed10,Tawarmalani04}).
An OBBT procedure receives in input a box $B$ and returns a tightened box in output,
removing feasible points which do not allow to improve the current best feasible solution. More formally, let ${\cal B}$ be the set of $n$-dimensional boxes. Then:
$$
OBBT:{\cal B} \rightarrow {\cal B}\ : \ OBBT(B)\subseteq B\ \ \ \mbox{and}\ \ \ F({\bf x},{\bf y}, v)\geq GUB\ \ \ \ \forall\ ({\bf x},{\bf y})\in [B\cap P_{{\bf x},{\bf y}}] \setminus OBBT(B).
$$
In our approach, we propose to employ an OBBT procedure, which is expensive but, as we will see, also able to considerably reduce the number of branch-and-bound nodes.
The lower and upper limits $\ell_{x_i}, \ell_{y_i},  u_{y_i},  u_{x_i}$, $i=1,\ldots,n$ are refined through the solution of LP problems having the feasible set defined by constraints (\ref{simp})-(\ref{MC2}) and the additional constraint
\begin{equation}
\label{eq:optbasred}
\frac{1}{2}\sum_{i=1}^n f_i + v  \le GUB,
\end{equation}
stating that we are only interested at feasible solutions where the underestimating function, i.e., the left-hand side of the constraint, corresponding to the objective function (\ref{obj}), is not larger than the current upper bound $GUB$. 
Thus, each call of this OBBT procedure requires the solution of $4n$ LP problems with the following objective functions:
$$
\ell_{x_i}/ u_{x_i}=\min/\max\ x_i,\ \ \ \ell_{y_i}/ u_{y_i}=\min/\max\ y_i,\ \ \ i=1,\ldots,n.
$$
Note that all these problems are bounded in view of the fact
that ${\bf x}$ is constrained to belong to the unit simplex.
In fact, what we observed through our computational experiments is that it is not necessary to solve all $4n$ LPs
but it is enough to concentrate the effort on the most 'critical' variables. More precisely, in order to reduce the number of LPs without compromising the performance, we employed the following strategies (see also \cite{Gleixner17} for strategies to reduce the effort).
Taking into account the quantities $g_i$ computed in (\ref{gap_i}), we notice that the larger the $g_i$ value, the higher is the need for a more accurate underestimation of the corresponding bilinear term.
Then, we solved the following LP problems.
\begin{itemize}
\item $\lceil 0.2 n \rceil$ LP problems with objective function $\min \ y_i$, for all $i$ corresponding to the $\lceil 0.2 n \rceil$ largest $g_i$ values;
\item a fixed number $\lceil 0.1 n \rceil$ of LP problems with
objective function $\max \ y_i$, for all $i$ corresponding to the $\lceil 0.1 n \rceil$ largest $g_i$ values;
\item again $\lceil 0.1 n \rceil$ LP problems with
objective function $\max \ x_i$, for all $i$ corresponding to the $\lceil 0.1 n \rceil$ largest $g_i$ values;
\item no LP problem with objective function $\min\ x_i$.
\end{itemize}
These choices have been driven by some experimental observations. In particular, we noticed that the lower limit for $y_i$ is the most critical
for the bound computation or, stated in another way, constraint
$$
 f_i \geq \ell_{y_i}x_i + \ell_{x_i} y_i -\ell_{x_i} \ell_{y_i},
$$
is often the active one. For this reason  a larger budget of LP problems is allowed to improve this lower limit with respect to the upper limits. Instead, we never try to improve the lower limit $\ell_{x_i}$ because it is experimentally observed that
this limit can seldom be improved.\newline
This way, the overall number of LPs to be solved at each call of the OBBT procedure is reduced to approximately $0.4n$.
Note that rather than solving all LP problems with the same feasible set, we could solve each of them with a different feasible region by incorporating all previously computed new limits in the definition of the feasible region for the next limit to be computed. That leads to sharper bounds, however we excluded this opportunity since we observed that LP solvers strongly benefit from the opportunity
of solving problems over the same feasible region. 
\newline\newline\noindent
The underestimating function depends on the lower and upper limits $ \ell_{x_i}, \ell_{y_i},  u_{y_i},  u_{x_i}$.
Thus, once we have updated all such limits, we can call procedure OBBT again in order to further reduce the limits. These can be iteratively reduced until some stopping condition is fulfilled. Such iterative procedure has been proposed and theoretically investigated, e.g., in \cite{DomRed10}.
That obviously increases the computational cost per node, since
the overall number of LPs to be solved at each node is now approximately $0.4n$ times the number of calls to the procedure OBBT, which depends on the stopping condition. But, again, we observed that the additional computational cost is compensated by a further reduction of the overall number of nodes in the branch-and-bound tree.
\newline\newline\noindent
It is important to stress at this point that OBBT procedures in general and the one proposed here in particular, are not new in the literature. The main contribution of this work lies in the observation that a very aggressive application of the proposed OBBT, while increasing considerably the computational cost per node, is the real key for an efficient solution of the addressed problem. Indeed, we will see through the experiments, that our approach is able to significantly outperform commercial QP solvers like {\tt CPLEX} and {\tt GUROBI}, and a solver like {\tt BARON}, which is strongly based on tightening techniques. This fact suggests that the intensive application of OBBT procedures might enhance the performance  of QP solvers not only over the addressed QP problems but also over more general ones.
\subsection{Pseudo-code of the branch-and-bound approach}
In this section we collect all the previously described tools and present the pseudo-code of the proposed branch-and-bound approach. In Line 1 an initial box $B_0$ is introduced and the collection of branch-and-bound nodes ${\cal C}$ still to be explored is initialized with it. In Line 2 a lower bound over $B_0$ is computed, while in Line 3 the current best observed feasible point ${\bf z}^\star$ and the current $GUB$ value are initialized. Take into account that such values can also be initialized after running a few local searches from randomly generated starting points.
Lines 4--21 contain the main loop of the algorithm. Until the set of nodes still to be explored is not empty, the following operations are performed. In Line 5 one node in ${\cal C}$ with the lowest lower bound is selected. In Line 6 the index $k$ of the branching variable is selected as the one with the largest gap  $g_i$ as defined in (\ref{gap_i}). In Line 7 the branching operation is performed. In Line 8 the selected node $\bar{B}$ is removed from ${\cal C}$, while in Lines 9--19 the following operations are performed for each of its child nodes. In the loop at Lines 10--17, first procedure OBBT is applied and then the lower bound over the tightened region is computed, until a stopping condition is satisfied. In particular, in our experiments we iterate until the difference between the (non-decreasing) lower bounds at two consecutive iterations fall below a given threshold $\epsilon$ ($\epsilon =10^{-3}$ in our experiments). In Lines 13--16 both ${\bf z}^\star$ and $GUB$ are possibly updated through the optimal solution of the relaxed problem.
In Line 18 we add the child node to ${\cal C}$. Finally, in Line 20 we remove from ${\cal C}$ all nodes with a lower bound
not lower than  $ (1-\varepsilon)GUB$, where $\varepsilon$ is a given tolerance value.
In all the experiments, we fixed a relative tolerance $\varepsilon=10^{-3}$, which is considered adequate for practical applications.
\newline\newline\noindent
Note that we do not discuss convergence of the proposed branch-and-bound approach, since it easily follows by rather standard and general arguments which can be found in \cite{HorstTuy93}.
\newline\newline\noindent
In Algorithm \ref{Alg:BandB} we highlighted with a frame box, both the stopping condition in Line 10 and the call to the OBBT procedure at Line 11, since the performance of the proposed algorithm mainly depends on how these two lines are implemented.

In what follows, in order to stress the importance of bound-tightening, we will refer to the proposed approach as Branch-and-Tight (B\&T), which belongs to the class of Branch-and-Cut approaches, since tightening the bound of a variable is a special case of cutting plane. 
\begin{algorithm}
\caption{Branch-and-bound algorithm}
\label{Alg:BandB}
\AlgData{${\bf S}\in \mathbb{R}^{n\times n}$, ${\bf A}\in \mathbb{R}^{m\times n}$,  $\varepsilon>0$\;}
Let $B_0$ be an initial box and set ${\cal C}=\{B_0\}$  \label{l1}\;
Compute the lower bound $L(B_0)$ through (\ref{boundLP})  \label{l2} \; 
${\bf z}^\star={\bf x}^\star(B_0)$ and $GUB=F\left({\bf x}^\star(B_0),{\bf y}^\star(B_0),v^\star(B_0)\right)$  \label{l3}\;
\While{${\cal C}\neq \emptyset$  \label{l5}}{
$\bar{B}\in \arg\min_{B\in {\cal C}} L(B)$  \label{l6}\;
$k\in \arg\max_{i=1,\ldots,n} x_i^\star(\bar{B}) y_i^\star(\bar{B})-f_i^\star(\bar{B})$  \label{l7}\;\vspace{0.05cm}
Branch $\bar{B}$ into $B_1=\bar{B}\cap \{y_k\leq  y_k^\star(\bar{B})\}$ and $B_2=\bar{B}\cap \{y_k\geq  y_k^\star(\bar{B})\}$ \label{l8}\;
${\cal C}={\cal C}\setminus \{\bar{B}\}$ \label{l9}\; 
\For{$i\in \{1,2\}$  \label{l10}}{
\While{\framebox{A stopping condition is not satisfied}  \label{l11}}{
\framebox{$B_i=OBBT(B_i)$} \label{l12}\;
Compute the lower bound $L(B_i)$ through (\ref{boundLP})  \label{l13}\;
\If{$F({\bf x}^\star(B_i), {\bf y}^\star(B_i),v^\star(B_i))<GUB$ \label{l14}}{
${\bf z}^\star={\bf x}^\star(B_i)$  \label{l15}\;
$GUB=F({\bf x}^\star(B_i), {\bf y}^\star(B_i),v^\star(B_i))$  \label{l16}\;
\label{l17}}
 \label{l18}}
${\cal C}={\cal C} \cup \{B_i\}$ \label{l20}\;
 \label{l22}}
${\cal C}={\cal C}\setminus  \{B\in {\cal C}\ :\ L(B)\geq (1-\varepsilon)GUB\}$  \label{l23}\;
}
\Return{${\bf z}^\star$, $GUB$\;}
\end{algorithm}
$\ $\newline\newline\noindent

\section{Numerical Results}
\label{sec:numrel}
\subsection{Test instances description}
\label{sec:testinst}
The game is about spot checking a set of $n$ places to guard them against an adversary. The places are spatially scattered, with a directed weighted graph describing the connections (direct reachability) of place $v$ from place $u$ by an edge $v\to u$ with a random length.
\newline\newline\noindent
The payoffs in the game are given by an $n\times n$ matrix ${\bf A}$ (so $m=n$ in the above description), and are interpreted as the loss that the defending player 1 suffers when checking place $i$ while the attacker is at place $j$.
Thus, the defender can:
\begin{itemize}
    \item either miss the attacker ($i\neq j$) in which case there will be a Weibull-distributed random loss with shape parameter 5 and scale parameter 10.63 (so that the variance is 5);

    \item or hit the attacker at $i=j$,  in which case there is zero loss.
\end{itemize}
The defender is thus minimizing, and the attacker is maximizing. The problem above is that of the defender. The Nash equilibrium then gives the optimal random choice of spot checks to minimize the average loss.
To avoid trivialities, the payoff matrices are constructed not to admit pure strategy equilibria, so that the optimum (without switching cost) is necessarily a mixed strategy.
\newline\newline\noindent
As for the switching cost, if the defender is currently at position $i$ and next -- according to the optimal random choice -- needs to check the (non-adjacent) place $j$, then the cost for the switch from $i$ to $j$ is the shortest path in the aforementioned graph (note that, since the graph is directed, the matrix ${\bf S}$ is generally nonsymmetric).
\newline\newline\noindent
For the random instances, the matrix ${\bf S}$ is thus obtained from a (conventional) all-shortest path algorithm applied to the graph. Note that the graph is an Erd\"{o}s-Renyi type graph with $n$ nodes and $p=0.3$ chance of any two nodes having a connection.

\begin{remark} The Erd\"{o}s-Renyi model is here a suitable  description of patrolling situations in areas where moving from any point to any other point is possible without significant physical obstacles in between. Examples are water areas (e.g., coasts) or natural habitats (woods, \ldots), in which guards are patrolling. It goes without saying that implementing the physical circumstances into the patrolling problem amounts to either a particular fixed graph topology or class of graphs (e.g., trees as models for harbor areas, or general scale-free networks describing communication relations). Such constrained topologies, will generally induce likewise constrained and hence different (smaller) strategy spaces, but leave the problem structure as such unchanged, except for the values involved.
\end{remark}
The weights in the graph were chosen exponentially distributed with rate parameter $\lambda=0.2$, and the Weibull distribution for the losses has a shape parameter 5 and scale parameter $\sim 10.63$, so that both distributions have the same variance of 5.
\begin{remark} 
The choice of the Weibull distribution is because of its heavy tails, useful to model extreme events (in actuarial science, where it appears as a special case of the generalized extreme value distribution). If the graph is an attack graph, we can think of possibly large losses that accumulate as the adversary traverses an attack path therein (but not necessarily stochastically independent, which the Weibull-distribution sort of captures due to its memory property). Besides, both the exponential and the Weibull distribution only take non-negative values, and thus lend themselves to a meaningful assignment of weights as ``distances'' in a graph. 
\end{remark}
The graph sizes considered are $n=50,75,100$ and for each graph size we consider ten random instances.
We restricted the attention to $\alpha$ values in $\{0.3,0.4,\ldots,0.9\}$ since problems with $\alpha$ values smaller than 0.3 and larger than 0.9 turned out to be simple ones. The overall number of instances is, thus, 210 (70 for each size $n=50,75,100$).
\newline\newline\noindent
Note that all the data of the test instances are available at the web site \url{http://www.iasi.cnr.it/~liuzzi/StQP}.

\subsection{Description of the experiments}
\label{sec:expliter}
The problem discussed in this paper belongs to the class
of nonconvex QP problems with linear constraints, which is a quite active research area.
Even well-known commercial solvers, like {\tt CPLEX} and {\tt GUROBI}, have recently included the opportunity of solving these problems. 
In \cite{Xia19} different solution approaches have been compared over different nonconvex QP problems, namely: Standard Quadratic Programming problems (StQP), where the feasible region is the unit simplex; BoxQP, where the feasible region is a box; and general QPs, where the feasible region is a general polytope (in \cite{Bonami18} an extensive comparison has also been performed more focused on BoxQPs). 
The approaches tested in \cite{Xia19} have been the nonconvex QP solver of {\tt CPLEX}, {\tt quadprogBB} (see \cite{ChenBurer12}), {\tt BARON} (see \cite{Sahinidis96}), and {\tt quadprogIP}, introduced in \cite{Xia19}. According to the computational results reported in that work, solvers {\tt quadprogIP} and {\tt quadprogBB} have quite good performance on some subclasses. More precisely, {\tt quadprogIP} works well on StQP problem (see also \cite{Gondzio18} for another approach working well on this subclass), while {\tt quadprogBB} performs quite well on BoxQPs, especially when the Hessian matrix of the objective function is dense. However, they do not perform very well on QP problems with general linear constraints. Some experiments we performed show that their performance is not good also on the QP subclass discussed in this paper. For this reason we do not include their results in our comparison. Thus, in the comparison we included: the nonconvex QP solver of {\tt CPLEX}, the best performing over QPs with general linear constraints according to what reported in \cite{Xia19}; the nonconvex QP solver of {\tt GUROBI}, which has been recently introduced and is not tested in that paper; {\tt BARON}, since bound-tightening, which, as we will see, is the most relevant operation in the proposed approach, also plays a central role in that solver.
\newline\newline\noindent
We performed three different sets of experiments: 
\begin{itemize}
    \item Experiments to compare our approach B\&T with the above mentioned existing solvers over the subclass of QP problems discussed in this paper (only at dimension $n=50$, which, as we will see, is already challenging for all the competitors);
    \item Experiments with B\&T by varying the intensity of bound-tightening (no bound-tightening, light bound-tightening, strong bound-tightening), in order to put in evidence that (strong) bound-tightening is the key operation in our approach;
    \item Experiments with B\&T at dimensions $n=50,75,100$, in order to see how it scales as the dimension increases.
\end{itemize}
 All the experiments have been carried out on an Intel\textsuperscript{\textregistered} Xeon\textsuperscript{\textregistered} gold 6136 CPU at 3GHz with 48 cores and 256GB main memory. 
 The algorithm has been coded using the Julia \cite{julia_paper} language (version 1.3.1). Doing the implementation we parallelized as much as possible the 
bound-tightening procedure discussed in Section \ref{sec:boundtight}, where many LPs with the same feasible region need to be solved. The code is available for download at the URL \url{http://www.iasi.cnr.it/~liuzzi/StQP}.
\subsubsection{Comparison with the existing literature}
 As a first experiment, we compare our method with the commercial solvers {\tt BARON}, {\tt CPLEX} and {\tt Gurobi}. We run all these methods over ten instances at dimension $n=50$ with $\alpha\in \{0.3,0.4,\ldots,0.9\}$ (thus, overall, 70 instances). We set a time limit of 600 seconds. A relative tolerance $\varepsilon=10^{-3}$ is required for all solvers.
 In Table \ref{tab:comsolv}, we report the average performance. For each method we report the number of nodes (column {\tt nn}), the percentage gap after the time limit and in brackets the computational time needed to reach it (column {\tt 	GAP \% (s)}), and finally the percentage of success, i.e. the percentage of instances solved to optimality within the time limit (column {\tt Succ \%}). In our opinion, this table reports the most important finding of this paper. It can be seen that all the commercial solvers fail on most of the instances (apart from a small number with $\alpha=0.9$), whereas our approach solves all the instances with an average time of less than 30 seconds (the complete results are reported in  Appendix \ref{sec:detres}).
 These results show that, although commercial solvers are fully developed, there is still room for improvements. In particular, it seems that performing bound-tightening in a very intensive way can  strongly enhance the performance of a solution approach. It might be the case that even commercial solvers may strongly benefit from an intensive application of this procedure. In fact, as previously recalled, {\tt BARON} already incorporates bound-tightening techniques but, as we will see in the following set of experiments, the intensity with which bound-tightening is applied also makes a considerable difference.  
 
 {\tiny
\begin{table}[ht]
    \centering
    \begin{tabular}{|l|l|l|l|l|l|l|l|l|l|l|l|l|}\hline
$\alpha$	&	\multicolumn{3}{|c|}{BARON}		&\multicolumn{3}{|c|}{CPLEX}					&	\multicolumn{3}{|c|}{GUROBI}&	\multicolumn{3}{|c|}{B\&T}\\\hline
	&	nn	&	GAP \% (s)	&	\% 	&	nn	&	GAP \% (s)	&	\% 	&	nn	&	GAP \% (s)	&	Succ \% 	&	nn	&	GAP \% (s)	&	Succ \% 	\\\hline
	0.3	&	555	&	4.27(600)	&	0	&	254487	&	1.59(600)	&	0	&	1072784	&	3.01(600)	&	0	&	11.2	&	0 (6.701)	&	100	\\\hline
	0.4	&	511	&	6.33(600)	&	0	&	241266	&	2.79(600)	&	0	&	933319	&	5.03(600)	&	0	&	70	&	0 (32.576)	&	100	\\\hline
	0.5	&	496	&	8.91(600)	&	0	&	233058	&	4.81(600)	&	0	&	912301	&	7.25(600)	&	0	&	80.6	&	0 (32.661)	&	100	\\\hline
	0.6	&	432	&	11.29(600)	&	0	&	936364	&	7.37(600)	&	0	&	990219	&	9.59(600)	&	0	&	45.8	&	0 (18.46)	&	100	\\\hline
	0.7	&	449	&	14.2(600)	&	0	&	1040716	&	9.73(600)	&	0	&	1118815	&	11.93(600)	&	0	&	57.6	&	0 (17.973)	&	100	\\\hline
	0.8	&	350	&	16.31(600)	&	0	&	1060996	&	10.51(600)	&	0	&	1323985	&	11.94(600)	&	0	&	54.6	&	0 (11.199)	&	100	\\\hline
	0.9	&	156	&	4.79(491.65)	&	40	&	1012193	&	2.31(444.41)	&	50	&	650214	&	0.27(600)	&	90	&	61.6	&	0 (7.562)	&	100	\\\hline
\end{tabular}
\caption{Average performance of all the solvers on ten instances for each value of $\alpha$ when $n=50$. The column $nn$ represents the average number of nodes, the column 	$GAP \% (s)$ reports the percentage GAP after at most 600 seconds and in brackets the average CPU time in seconds. The column $Succ \%$ represents the percentage of success among the ten instances.}\label{tab:comsolv}
\end{table}}

\subsubsection{Importance of bound-tightening}
As already stressed many times, the quite good performance of B\&T is due to the bound-tightening procedure. It is now time to show it with numbers. To this end, besides the already proposed setting for our approach, we ran it under two different settings: 
\begin{description}
\item[{\bf No bound-tightening}] at each node we do not apply the procedure $OBBT$, but we simply compute the lower bound by solving problem \eqref{boundLP};
\item[{\bf Light Bound-Tightening}] at each node we only solve the following LPs once (and, consequently, we compute the solution of problem \eqref{boundLP} only once)  \begin{itemize}
\item $\lceil 0.1 n \rceil$ LP problems with objective function $\min \ y_i$, for all $i$ corresponding to the $\lceil 0.1 n \rceil$ largest $g_i$ values;
\item a fixed number $\lceil 0.05 n \rceil$ of LP problems with
objective function $\max \ y_i$, for all $i$ corresponding to the $\lceil 0.05 n \rceil$ largest $g_i$ values;
\item again $\lceil 0.05 n \rceil$ LP problems with
objective function $\max \ x_i$, for all $i$ corresponding to the $\lceil 0.05 n \rceil$ largest $g_i$ values;
\item no LP problem with objective function $\min\ x_i$.
\end{itemize}
\end{description}
Of course, this strongly reduces the effort per node. With no bound-tightening a single LP is solved per node, while with light bound-tightening $\lceil 0.2 n \rceil+1$ LPs are solved at each node. In fact, light bound-tightening already requires a considerable computational effort per node (and, as we will see, it is already enough to perform better than existing solvers). However, the originally proposed strong bound-tightening procedure, where the OBBT procedure is iteratively applied and at each iteration $\lceil 0.4 n \rceil$ LPs are solved, delivers better results.
In Table \ref{tab:bbt}, we report the average performance on the instances with $n=50$ in terms of number of nodes, number of LPs solved, CPU time in seconds and percentage gap of the three levels of bound-tightening. It is evident from the table that the OBBT procedure is what really makes the difference: most of the instances are not solved without bound-tightening, whereas the number of nodes and the CPU time needed to solve the instances decrease as we increase the level of bound-tightening. In Appendix \ref{sec:detres} we also report the full table with all the results.

{\tiny
\begin{table}[ht]
    \centering
    \begin{tabular}{|l|l|l|l|l|l|l|l|l|l|l|l|l|}\hline
	&	\multicolumn{4}{|c|}{No Bound Tightening}							&	\multicolumn{4}{|c|}{Light Bound Tightening}&	\multicolumn{4}{|c|}{Strong Bound Tightening}\\\hline
$\alpha$	&	nn	&	\#LPs	&	time	&	GAP\%	&	nn	&	\#LPs	&	time	&	GAP\%	&	nn	&	\#LPs	&	time	&	GAP\%	\\\hline
0.3	&	106934.8	&	106937.8	&	600	&	1.899	&	423	&	19252.7	&	26.24	&	0	&	11.2	&	4665	&	6.7	&	0	\\\hline
0.4	&	111981.2	&	111984.2	&	600	&	2.93	&	2883.4	&	129670.2	&	169.431	& 1.3$\times 10^{-2}$	&	70	&	27077.5	&	32.58	&	0	\\\hline
0.5	&	117565.6	&	117568.6	&	600	&	4.19	&	1949.8	&	83921.9	&	108.489	&	0	&	80.6	&	26668.6	&	32.661	&	0	\\\hline
0.6	&	117942.4	&	117945.4	&	600	&	4.932	&	1106.8	&	45581.3	&	57.73	&	0	&	45.8	&	14882.5	&	18.46	&	0	\\\hline
0.7	&	113859.4	&	113862.4	&	600	&	6.601	&	813.8	&	30621.5	&	39.447	&	0	&	57.6	&	14604.9	&	17.973	&	0	\\\hline
0.8	&	90646.8	&	90649.8	&	600	&	6.77	&	307	&	10600.6	&	13.23	&	0	&	54.6	&	9247	&	11.199	&	0	\\\hline
0.9	&	89502.6	&	89505.6	&	531.13	&	1.48	&	183.2	&	5810.2	&	7.43	&	0	&	61.6	&	6477.2	&	7.56	&	0	\\\hline
\end{tabular}
\caption{Average performance on all the 10 instances for each value of $\alpha$ when $n=50$ for the different levels of bound tightening. }\label{tab:bbt}
\end{table}}

\subsubsection{Computational results over the proposed test instances as $n$ increases}
\begin{figure}
\centering
  \begin{subfigure}[b]{0.45\textwidth}
    \includegraphics[width=\linewidth]{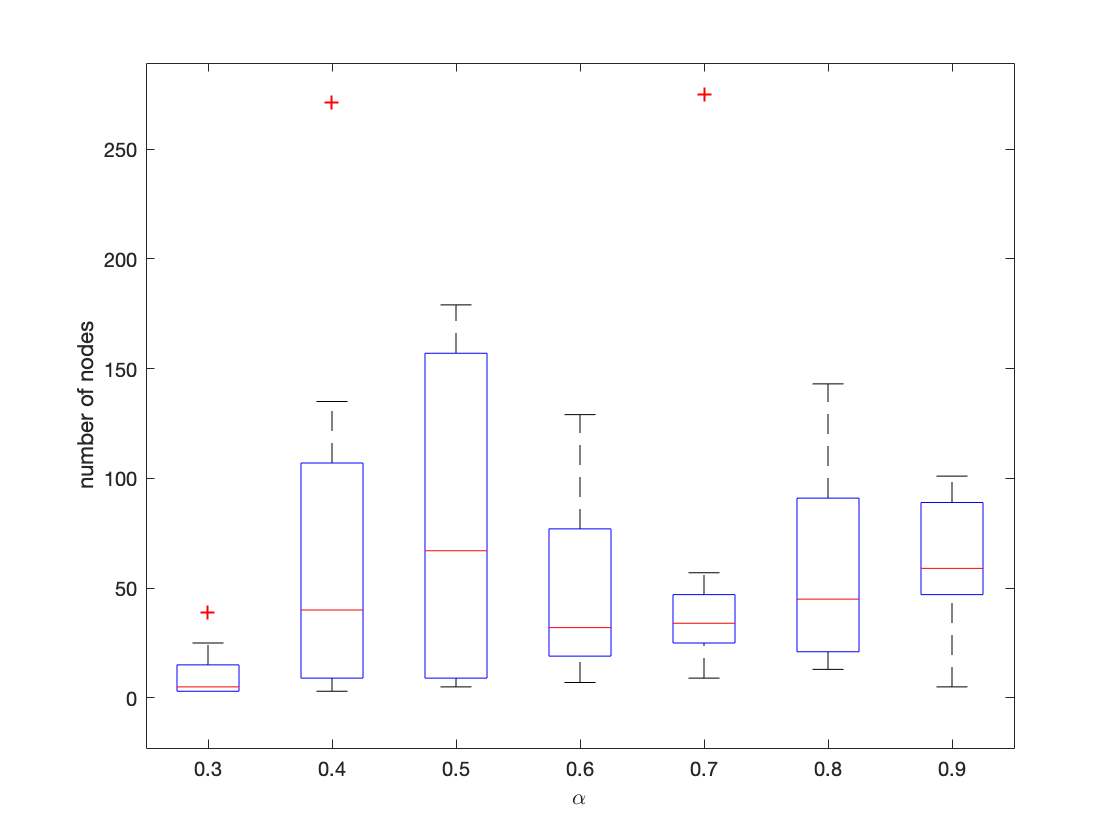}
    \caption{Total number of nodes for $n=50$}
    \label{fig:n50_nn}
  \end{subfigure}
  \begin{subfigure}[b]{0.45\textwidth}
  \includegraphics[width=\linewidth]{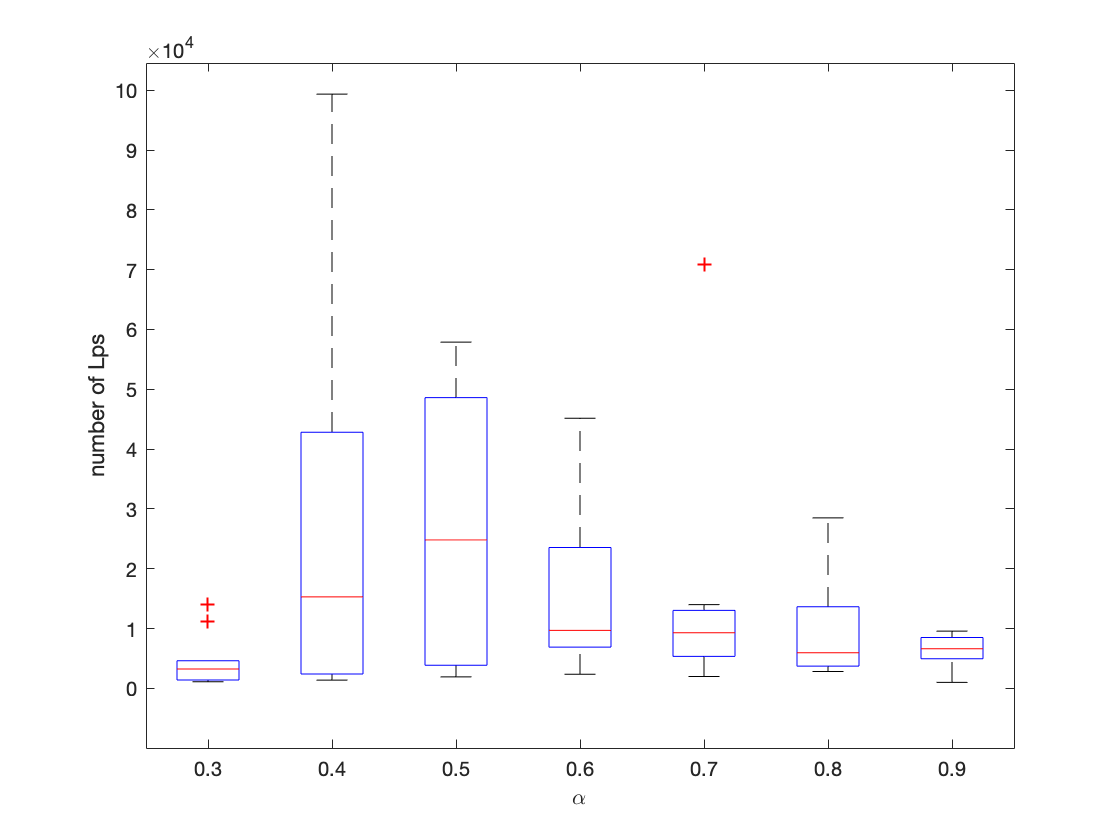}
    \caption{Number of LPs solved for $n=50$}
    \label{fig:n50_LP}
  \end{subfigure}
  \begin{subfigure}[b]{0.45\textwidth}
  \includegraphics[width=\linewidth]{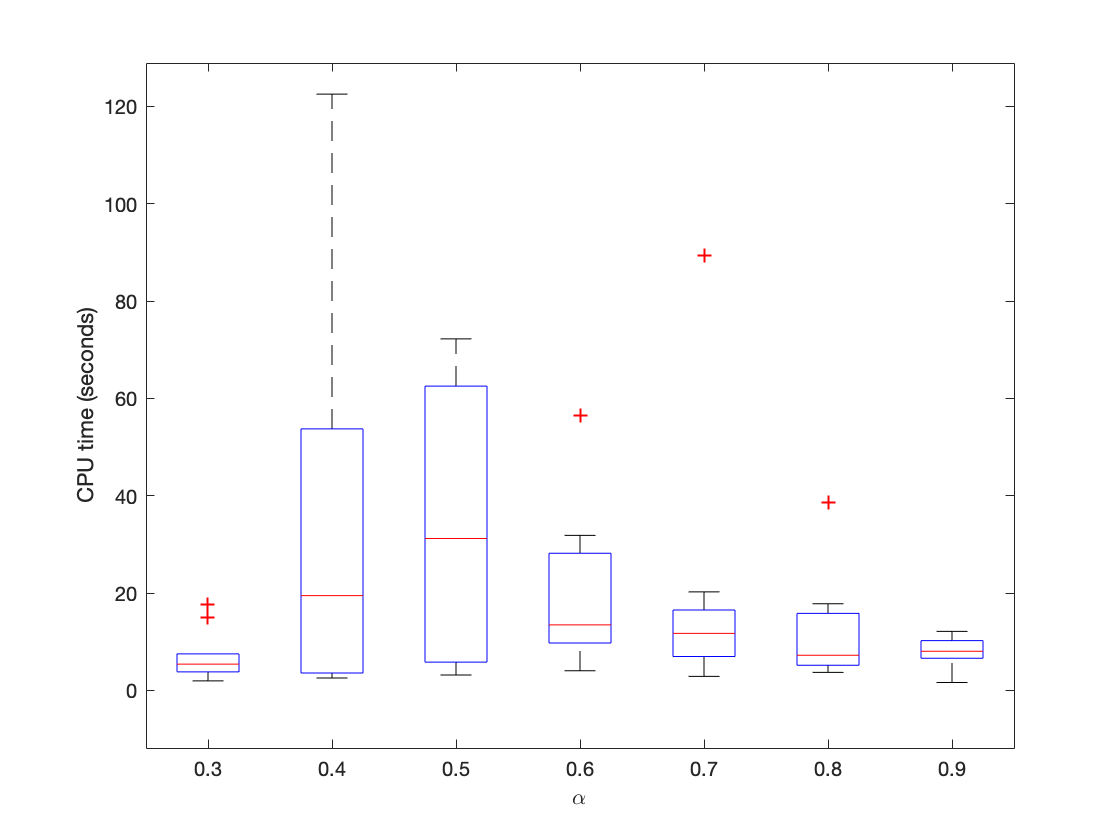}    \caption{CPU times in seconds}
    \label{fig:n50_sec}
  \end{subfigure}
   \caption{Box plots for different performance measures for $n=50$} \label{fig:n50}
\end{figure}

\begin{figure}
\centering
  \begin{subfigure}[b]{0.45\textwidth}
    \includegraphics[width=\linewidth]{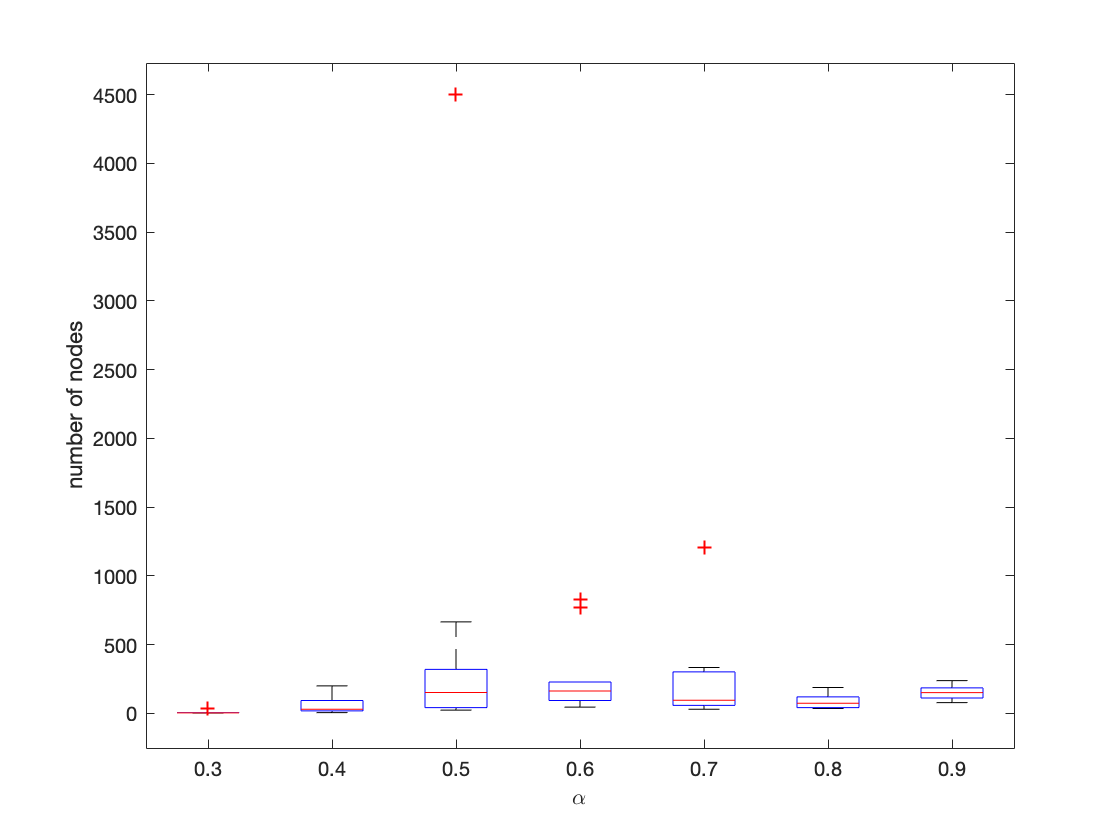}
    \caption{Total number of nodes for $n=75$  \label{fig:n75_nn}}
  \end{subfigure}
  \begin{subfigure}[b]{0.45\textwidth}
  \includegraphics[width=\linewidth]{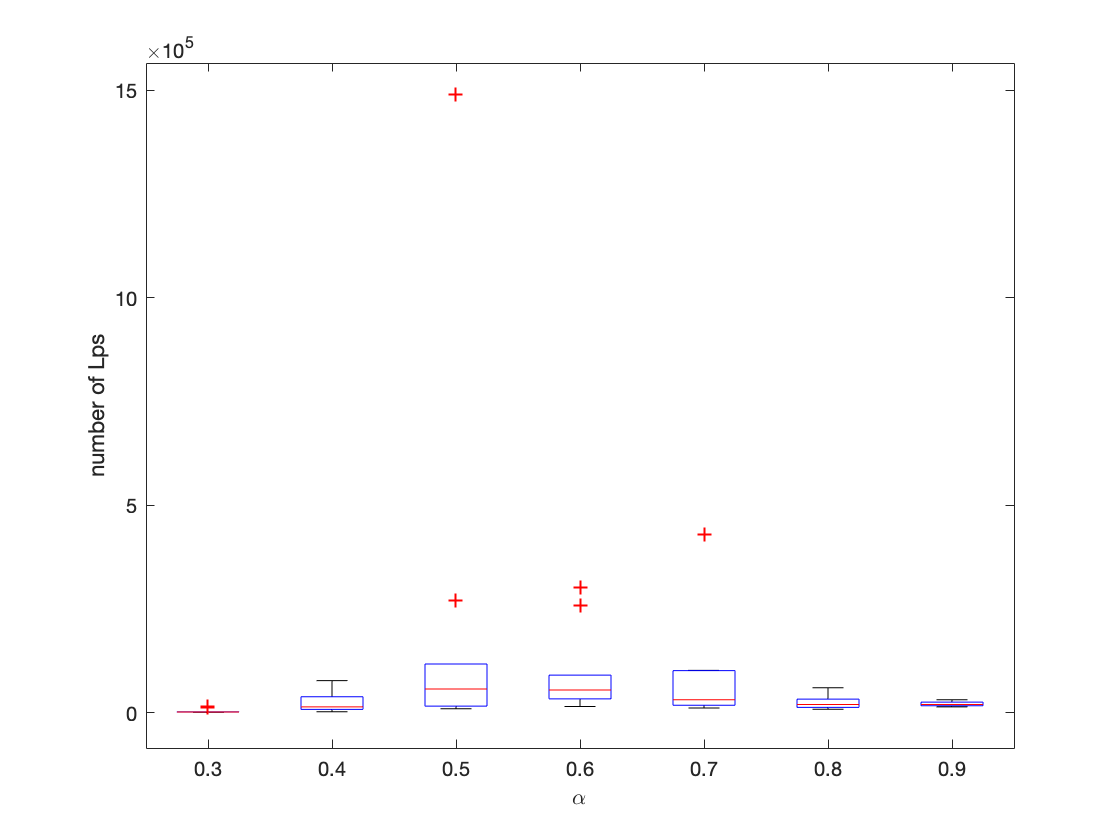}
    \caption{Number of LPs solved for $n=75$ \label{fig:n75_LP}}
  \end{subfigure}
  \begin{subfigure}[b]{0.45\textwidth}
  \includegraphics[width=\linewidth]{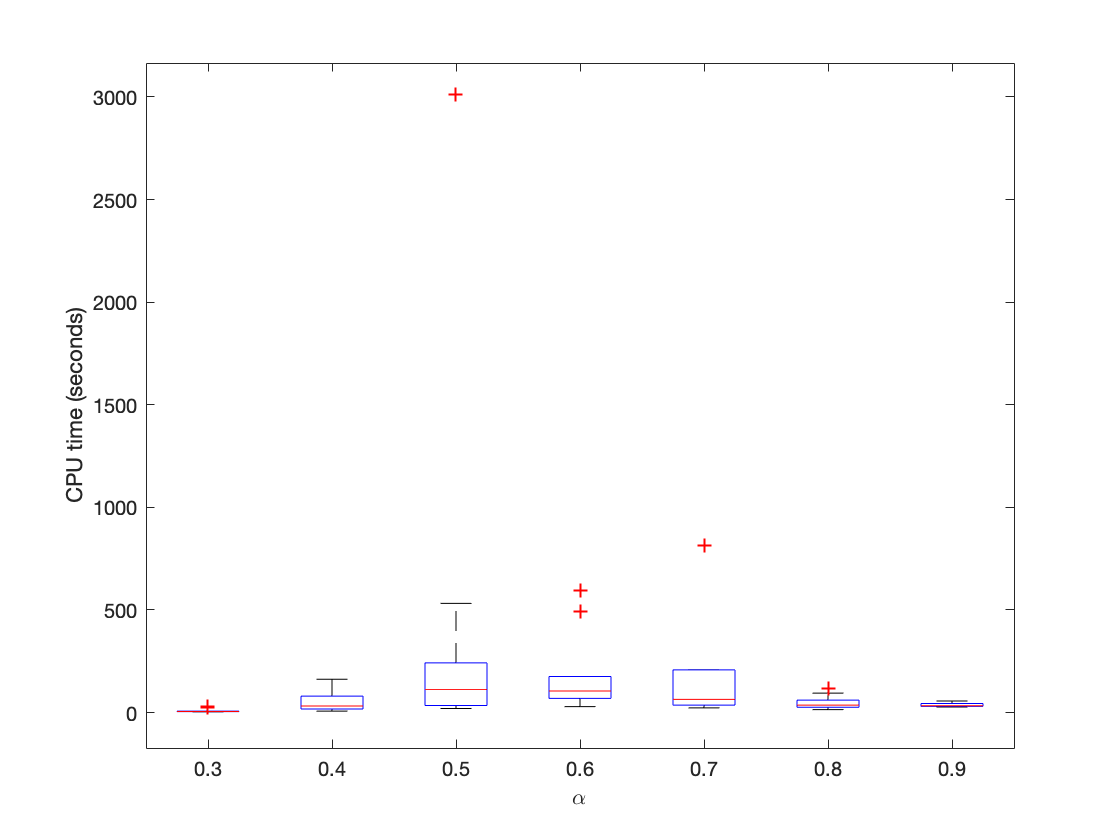}    \caption{CPU times in seconds \label{fig:n75cpu}}
  \end{subfigure}
   \caption{Box plots for different performance measures for $n=75$ \label{fig:n75} }
\end{figure}

\begin{figure}
\centering
  \begin{subfigure}[b]{0.45\textwidth}
    \includegraphics[width=\linewidth]{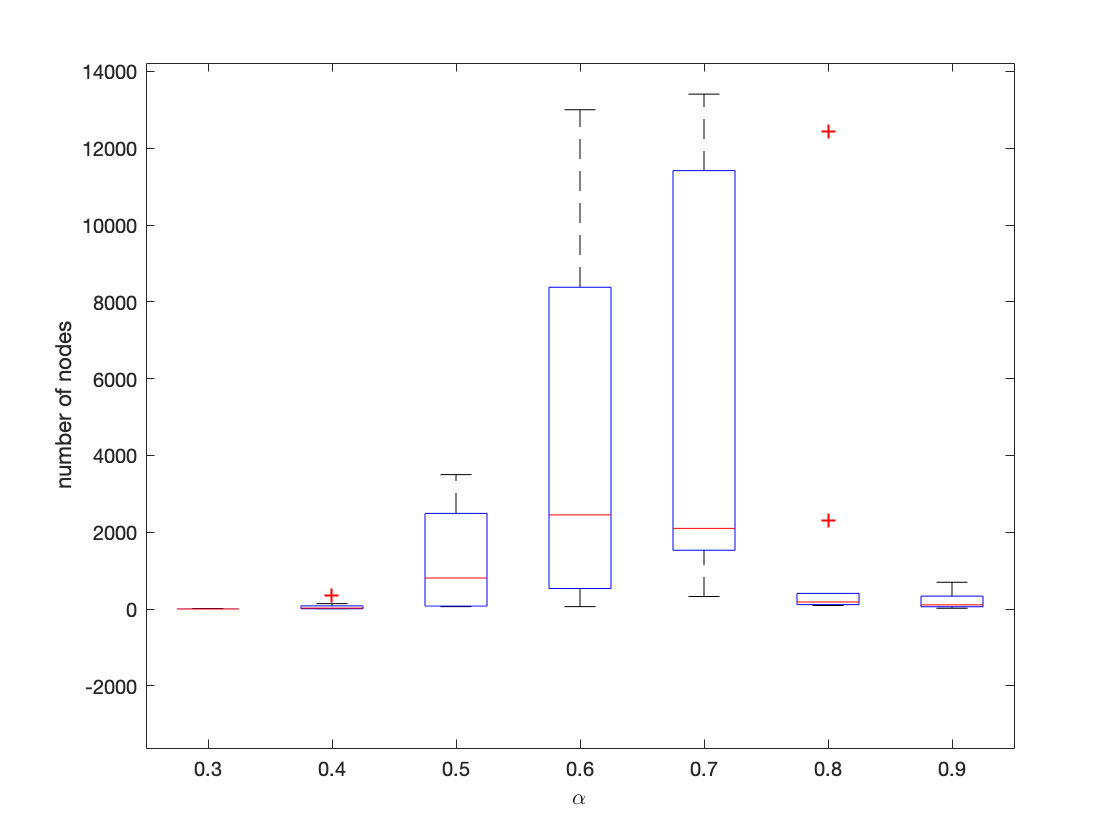}
    \caption{Total number of nodes for $n=100$}
    \label{fig:n100_nn}
  \end{subfigure}
  \begin{subfigure}[b]{0.45\textwidth}
  \includegraphics[width=\linewidth]{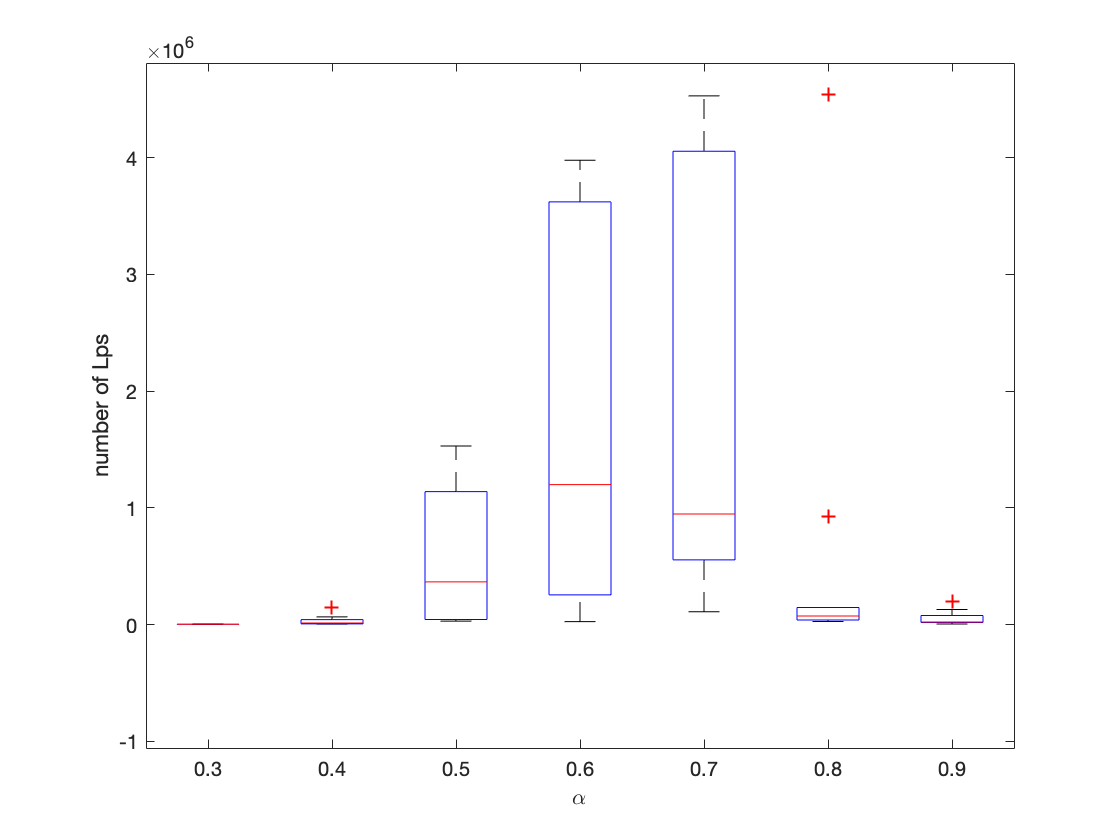}
    \caption{Number of LPs solved for $n=100$}
    \label{fig:n100_LP}
  \end{subfigure}
  \begin{subfigure}[b]{0.45\textwidth}
  \includegraphics[width=\linewidth]{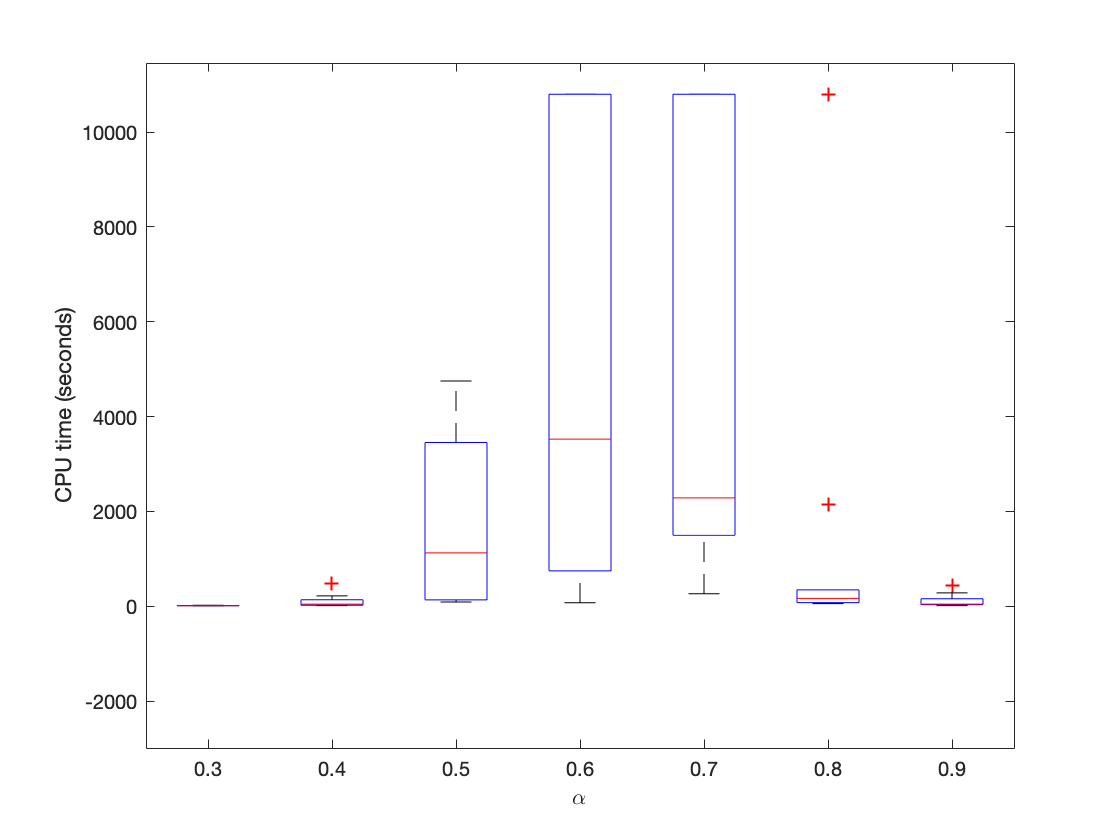}    \caption{CPU times in seconds}
    \label{fig:n100_sec}
  \end{subfigure}
  \begin{subfigure}[b]{0.45\textwidth}
  \includegraphics[width=\linewidth]{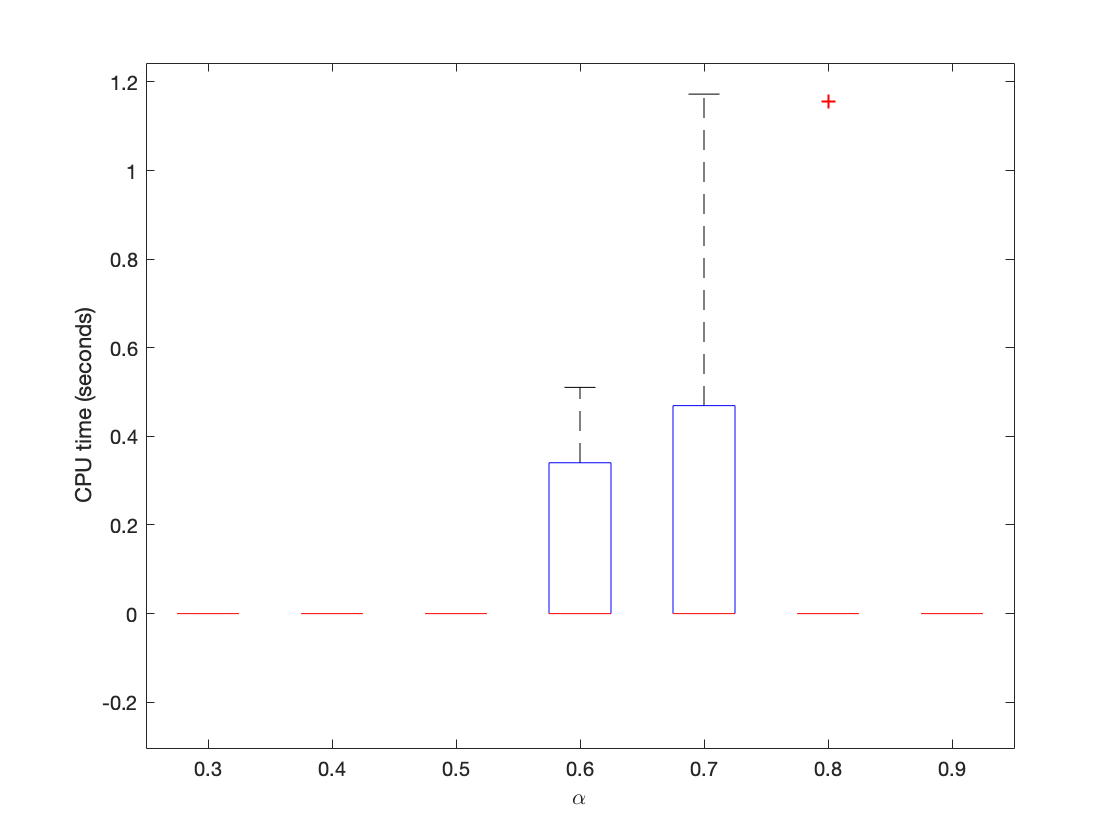}    \caption{\% Gap at the time limit}
    \label{fig:n100_gap}
  \end{subfigure}
   \caption{Box plots for different performance measures for $n=100$} \label{fig:n100}
\end{figure}

As a final experiment, we show the behavior of B\&T as the dimension $n$ increases. 
We have solved ten instances for three different sizes $n=50,~75,~100$ and the usual values of $\alpha=\{0.3,0.4,\ldots,0.9\}$ (thus, overall 210 instances). Note that, for all the different values of $n$, lower and larger values of $\alpha$ with respect to those we tested give rise to much simpler instances (recall that the problem becomes polynomially solvable for the extreme values $\alpha=0$ and $\alpha=1$).
We set a time limit of $10800$s for all instances. 
For $n=50$ and $n=75$ we solve all the instances to optimality (in fact, the largest time to solve an instance with $n=50$ is about 2 minutes, whereas for $n=75$ the largest time is below 1 hour, but most of the problems are solved within 10 minutes).  In Figures \ref{fig:n50_nn}-\ref{fig:n50_sec} we report the box plot of  number of nodes, number of LPs and CPU time needed for the different values of $\alpha$ when $n=50$. The figure shows that the hardest instances are the ones corresponding to the central values of $\alpha$ and this will turn out to hold true also at larger dimensions. We also observe that the overall number of nodes is extremely limited, thus confirming that, while computationally expensive, the bound-tightening procedure allows to considerably reduce the size of the branch-and-bound tree (again, this fact is observed also at larger dimensions).
In Figures \ref{fig:n75}, we report the different box plots for all the instances at $n=75$. It is worthwhile to remark that we solve most of them within ten minutes and exploring less than 300 nodes.
Finally, in Figure \ref{fig:n100} we report the performance of B\&T on problems of dimension $n=100$. In this case there are seven instances we are not able to solve within the time limit. These occur for values $\alpha\in\{0.6,0.7,0.8\}$, thus confirming that the central values of this parameter give rise to the most challenging instances. With respect to $n=50$ and $n=75$, we have the additional box plot displayed in Figure \ref{fig:n100_gap} reporting the final percentage gap when the time limit is reached. Note that it is never larger than 1.2\% and most of the times it is lower than 0.5\%, thus showing that, even when the algorithm does not terminate, the quality of the returned solution is guaranteed to be  high.

\section{Conclusions and future work}
\label{sec:concl}
In this paper we addressed some game theory problems arising in the context of network security. In these problems
there is an additional quadratic term, representing {\em switching costs}, i.e.. the costs for the defender of switching from a given strategy to another one at successive rounds of the game. The resulting problems can be reformulated as nonconvex QP with linear constraints.
Test instances of these problems turned out to be very challenging for existing solvers, and we propose to extend with them the current benchmark set of test instances for QP problems. 
We presented a spatial branch-and-bound approach to tackle these problems and we have shown that a rather aggressive application of an OBBT procedure is the key for their efficient solution.
The procedure is expensive, since it requires 
 multiple solutions of LP problems at each node of the branch-and-bound tree. But we empirically observed that 
 the high computational costs per node of the branch-and-bound tree are largely compensated by 
 the low number of nodes to be explored. As a topic for future research, we would like to further investigate the use of OBBT procedures in the solution of QP problems, and we would like to identify other cases, besides those addressed in this paper, where their intensive application may considerably enhance the performance of branch-and-bound approaches.

\bibliographystyle{spmpsci.bst}
\bibliography{secSTQP}

\begin{thebibliography}{10}
\providecommand{\url}[1]{{#1}}
\providecommand{\urlprefix}{URL }
\expandafter\ifx\csname urlstyle\endcsname\relax
  \providecommand{\doi}[1]{DOI~\discretionary{}{}{}#1}\else
  \providecommand{\doi}{DOI~\discretionary{}{}{}\begingroup
  \urlstyle{rm}\Url}\fi

\bibitem{Alpcan&Basar2010}
Alpcan, T., Ba{\c{s}}ar, T.: Network security: A decision and game-theoretic
  approach.
\newblock Cambridge University Press (2010)

\bibitem{Alpern.2011}
Alpern, S., Morton, A., Papadaki, K.: Patrolling games.
\newblock Operations Research \textbf{59}(5), 1246--1257 (2011)

\bibitem{julia_paper}
Bezanson, J., Edelman, A., Karpinski, S., Shah, V.B.: Julia: A fresh approach
  to numerical computing.
\newblock SIAM review \textbf{59}(1), 65--98 (2017).
\newblock \urlprefix\url{https://doi.org/10.1137/141000671}

\bibitem{Bonami18}
Bonami, P., G\"unl\"uk, O., Linderoth, J.: Globally solving nonconvex quadratic
  programming problems with box constraints via integer programming methods.
\newblock Mathematical Programming Computation \textbf{10}, 333--382 (2018)

\bibitem{DomRed10}
Caprara, A., Locatelli, M.: Global optimization problems and domain reduction
  strategies.
\newblock Mathematical Programming \textbf{125}(1), 123--137 (2010)

\bibitem{ChenBurer12}
Chen, J., Burer, S.: Globally solving nonconvex quadratic programming problems
  via completely positive programming.
\newblock Mathematical Programming Computation \textbf{4}(1), 33--52 (2012)

\bibitem{furini2019qplib}
Furini, F., Traversi, E., Belotti, P., Frangioni, A., Gleixner, A., Gould, N.,
  Liberti, L., Lodi, A., Misener, R., Mittelmann, H., et~al.: Qplib: a library
  of quadratic programming instances.
\newblock Mathematical Programming Computation \textbf{11}(2), 237--265 (2019)

\bibitem{Gleixner17}
Gleixner, A., Berthold, T., M\"uller, B., Weltge, S.: Three enhancements for
  optimization-based bound tightening.
\newblock Journal of Global Optimization \textbf{67}, 731--757 (2017)

\bibitem{Gondzio18}
Gondzio, J., Yildirim, E.A.: Global solutions of nonconvex standard quadratic
  programs via mixed integer linear programming reformulations.
\newblock arXiv preprint arXiv:1810.02307  (2018)

\bibitem{Hansen2012}
Hansen, K.A., Koucky, M., Lauritzen, N., Miltersen, P.B., Tsigaridas, E.P.:
  Exact algorithms for solving stochastic games.
\newblock In: Proceedings of the forty-third annual ACM symposium on Theory of
  computing, pp. 205--214 (2011)

\bibitem{HorstTuy93}
Horst, R., Tuy, H.: Global optimization: Deterministic approaches (2nd
  edition).
\newblock Springer Science \& Business Media (2013)

\bibitem{lozovanu_multiobjective_2005}
Lozovanu, D., Solomon, D., Zelikovsky, A.: Multiobjective games and determining
  pareto-nash equilibria.
\newblock Buletinul Academiei de {\c{S}}tiin{\c{t}}e a Republicii Moldova.
  Matematica \textbf{3}, 115--122 (2005)

\bibitem{McCormick76}
McCormick, G.P.: Computability of global solutions to factorable nonconvex
  programs: Part i—convex underestimating problems.
\newblock Mathematical Programming \textbf{10}(1), 147--175 (1976)

\bibitem{Pardalos91}
P.M., P., S.A., V.: Quadratic programming with one negative eigenvalue is
  {NP}-hard.
\newblock Journal of Global Optimization \textbf{1}(1), 15--22 (1991)

\bibitem{rass_physical_2017}
Rass, S., Alshawish, A., Abid, M.A., Schauer, S., Zhu, Q., De~Meer, H.:
  Physical intrusion games—optimizing surveillance by simulation and game
  theory.
\newblock IEEE Access \textbf{5}, 8394--8407 (2017)

\bibitem{rass_password_2018}
Rass, S., K{\"o}nig, S.: Password security as a game of entropies.
\newblock Entropy \textbf{20}(5), 312 (2018)

\bibitem{Rass.2017}
Rass, S., K{\"o}nig, S., Schauer, S.: On the cost of game playing: How to
  control the expenses in mixed strategies.
\newblock In: International Conference on Decision and Game Theory for
  Security, pp. 494--505. Springer (2017)

\bibitem{rass_numerical_2014}
Rass, S., Rainer, B.: Numerical computation of multi-goal security strategies.
\newblock In: International Conference on Decision and Game Theory for
  Security, pp. 118--133. Springer (2014)

\bibitem{Sahinidis96}
Sahinidis, N.V.: Baron: A general purpose global optimization software package.
\newblock Journal of Global Optimization \textbf{8}(2), 201--205 (1996)

\bibitem{Tambe2012}
Tambe, M.: Security and game theory: algorithms, deployed systems, lessons
  learned.
\newblock Cambridge University Press (2011)

\bibitem{Tawarmalani04}
Tawarmalani, M., Sahinidis, N.V.: Global optimization of mixed-integer
  nonlinear programs: A theoretical and computational study.
\newblock Mathematical Programming \textbf{99}(3), 563--591 (2004)

\bibitem{Motzkin65}
T.S., M., E.G., S.: Maxima for graphs and a new proof of a theorem of
  {T}ur\'an.
\newblock Canadian Journal of Mathematics \textbf{17}(4), 533--540 (1965)

\bibitem{Wachter.2018}
Wachter, J., Rass, S., K{\"o}nig, S.: Security from the adversary’s
  inertia--controlling convergence speed when playing mixed strategy
  equilibria.
\newblock Games \textbf{9}(3), 59 (2018)

\bibitem{Xia19}
Xia, W., Vera, J.C., Zuluaga, L.F.: Globally solving nonconvex quadratic
  programs via linear integer programming techniques.
\newblock INFORMS Journal on Computing \textbf{32}(1), 40--56 (2020)

\end{thebibliography}

\appendix
\section{Proof of Theorem \ref{theo:compl}}
\label{sec:complexity}
In this section we consider the complexity of problem (\ref{eqn:nlp}). Such problem is a nonconvex QP with linear constraints. NP-hardness of QP problems has been established for different special cases like, e.g., the already mentioned StQP problems (see \cite{Motzkin65}) and the Box QP problems (see, e.g., \cite{Pardalos91}). However, due to its special structure, none of the known complexity results can be applied to establish the NP-hardness of problem (\ref{eqn:nlp}).  Thus, in what follows we formally prove that its corresponding decision problem is NP-complete.
Let
\begin{equation}
\label{eq:objfun}
f({\bf x})=\frac{1}{2} {\bf x}^T {\bf S} {\bf x}+ \max_{k=1,\ldots,m} {\bf A}_k^T {\bf x},
\end{equation}
where $S_{ii}=0$ for each $i=1,\ldots,n$ and $S_{ij}\geq 0$ for each $i\neq j$, while ${\bf A}_k\geq {\bf 0}$, $k=1,\ldots,m$. Moreover, let
$$
\Delta_n=\{{\bf x}\in \mathbb{R}_+^n\ :\ {\bf e}^T{\bf x}=1\},
$$
be the $n$-dimensional unit simplex. After incorporating $\alpha$ and $(1-\alpha)$ respectively into ${\bf S}$ and ${\bf A}_k$, $k=1,\ldots,m$, problem (\ref{eqn:nlp})  is equivalent to $\min_{{\bf x}\in \Delta_n} f({\bf x})$. Then, we would like to establish the complexity of the following decision problem:
\begin{equation}
\label{eq:decprob}
\mbox{Given a constant $\xi\geq 0$}\ \ \ \exists\ {\bf x}\in  \Delta_n\ :\ f({\bf x})\leq \xi ?
\end{equation}
We prove the result by providing a polynomial transformation of the max clique decision problem: Given a graph $G=(V,E)$ and a positive integer $k\leq |V|$, does there exist a clique $C$ in $G$ with cardinality at least $k$?
We define the following instance of the decision problem (\ref{eq:decprob}).
Let
$$
S_{ij}=\left\{
\begin{array}{ll}
0 & \mbox{if $i=j$ or $(i,j)\in E$} \\[8pt]
n^4 & \mbox{otherwise.}
\end{array}
\right.
$$
Moreover, let $m=n$ and for each $k=1,\ldots,n$ let ${\bf A}_k={\bf e}_k$, where ${\bf e}_k$ is the vector with all components equal to 0, except the $k$-th one, which is equal to 1. Stated in another way, the piece-wise linear part
is $\max_{k=1,\ldots,n} x_k$. Finally, let $\xi=\frac{1}{k}$.
We claim that the minimum value of $f$ over $\Delta_n$ is not larger than $\xi=\frac{1}{k}$ if and only if $G$ contains a clique with cardinality $k$. The {\em if}  part is very simple. Indeed, let us consider the feasible solution $x_i=\frac{1}{k}$ if $i\in C$, where $C$ is a clique of cardinality $k$, and let $x_i=0$ otherwise. Then, the value of $f$ at this point is equal to $\frac{1}{k}$. Indeed , the value of the quadratic part is 0, while the value of the piece-wise linear part is $\frac{1}{k}$.
The proof of the {\em only if}  part is a bit more complicated. We would like to prove that, in case no clique with cardinality at least $k$ exists, then the minimum value of $f$ over the unit simplex is larger than $\frac{1}{k}$. Let us denote by ${\bf x}^*$ the minimum of $f$ over the unit simplex. Let
$$
K=supp({\bf x}^*)=\{i\ :\ x_i^*>0\},
$$
and let $C$ be the maximum clique over the sub-graph induced by $K$, whose cardinality is at most $k-1$.
We first remark that if, for some $i\in K\setminus C$, it holds that
$$
x_i^*\geq \frac{1}{n^2}\ \ \ \mbox{and}\ \ \ \sum_{j\in C\ :\ (i,j)\not\in E} x_j^*\geq \frac{1}{n^2},
$$
then the quadratic part contains the term
$$
n^4x_i^*\left( \sum_{j\in C\ :\ (i,j)\not\in E} x_j^*\right)\geq 1,
$$
which concludes the proof. Therefore, for $i\in K\setminus C$ we assume that either $x_i^*<\frac{1}{n^2}$ or
\begin{equation}
\label{eq:auxineq}
 \sum_{j\in C\ :\ (i,j)\not\in E} x_j^* < \frac{1}{n^2}.
\end{equation}
Now, let
$$
K_1=\left\{i\ :\ i\in K\setminus C\ \mbox{and}\ x_i^*\geq \frac{1}{n^2}\right\}.
$$
If $\exists\ k_1,k_2\in K_1$ and $(k_1,k_2)\not \in E$, then $n^4 x_{k_1}^*x_{k_2}^*\geq 1$, which concludes the proof.
Then, we assume that for each $k_1, k_2\in K_1$, $(k_1,k_2)\in E$, i.e., $K_1$ itself is a clique. Now let us consider the following subset of $C$
$$
C_1=\left\{i\in C\ :\ (i,k)\not\in E\ \mbox{for at least one $k\in K_1$}\right\}.
$$
It must hold that $|C_1|\geq |K_1|$. Indeed, if $|C_1|<|K_1|$, then $(C\setminus C_1) \cup K_1$ is also a clique with cardinality larger than $C$, which is not possible in view of the fact that $C$ has maximum cardinality.
Then, in view of (\ref{eq:auxineq}) we have that
$$
x_i^* < \frac{1}{n^2} \ \ \ \forall \ i\in C_1,
$$
and, moreover, by definition of $K_1$, we also have
$$
x_i^*< \frac{1}{n^2}\ \ \ \forall\ i\in K\setminus (K_1\cup C).
$$
Since $|C_1|\geq |K_1|$, we have that
$$
T=\left\{i\in K\ :\ x_i^*\geq \frac{1}{n^2}\right\},
$$
is such that $|T|\leq |K_1| + |C\setminus C_1|\leq |C_1|+|C\setminus C_1|=|C|\leq k-1$.
Thus, taking into account that
$$
\sum_{i\in K\setminus T} x_i^* < \frac{1}{n},
$$
we must have that
$$
\sum_{i \in T} x_i^* > 1-\frac{1}{n},
$$
and, consequently, taking into account that $|T|\leq k-1$, for at least one index $j\in T$ it must hold that
$$
x^*_j> \frac{1-\frac{1}{n}}{k-1}\geq \frac{1}{k},
$$
so that the piece-wise linear part of $f$ is larger than $\frac{1}{k}$, which concludes the proof.

\section{Detailed numerical results}\label{sec:detres}

{\tiny
\begin{table}[ht]
    \centering
    \begin{tabular}{|l|l|l|l|l|l|l|l|l|l|l|l|l|}\hline
$\alpha$	&	\multicolumn{3}{|c|}{BARON}		&\multicolumn{3}{|c|}{CPLEX}					&	\multicolumn{3}{|c|}{GUROBI}&	\multicolumn{3}{|c|}{Our approach}\\\hline
		&	nodes	&	time	&	GAP \%	&	 nodes	&	 time	&	GAP\% 	&	nodes 	&	time	&	GAP\%&	nodes 	&	time	&	GAP\%	\\\hline
0.3	&	503	&	600	&	3.94	&	264068	&	600	&	1.55	&	1088008	&	600	&	2.80	&	   3	&	8	&	0.0	\\\hline
0.3	&	514	&	600	&	3.05	&	306808	&	600	&	0.88	&	1325693	&	600	&	1.81	&	   5	&	4	&	0.0	\\\hline
0.3	&	625	&	600	&	5.29	&	190794	&	600	&	2.90	&	871168	&	600	&	4.16	&	   3	&	2	&	0.0	\\\hline
0.3	&	560	&	600	&	4.87	&	216867	&	600	&	2.08	&	987370	&	600	&	3.52	&	  25	&	14	&	0.0	\\\hline
0.3	&	550	&	600	&	4.53	&	260776	&	600	&	1.45	&	1060202	&	600	&	3.56	&	  39	&	17	&	0.0	\\\hline
0.3	&	525	&	600	&	5.14	&	269112	&	600	&	1.56	&	1071177	&	600	&	3.29	&	  15	&	6	&	0.0	\\\hline
0.3	&	584	&	600	&	4.57	&	197166	&	600	&	2.78	&	846598	&	600	&	3.50	&	   9	&	6	&	0.0	\\\hline
0.3	&	692	&	600	&	3.85	&	196197	&	600	&	1.65	&	1061279	&	600	&	2.79	&	   5	&	4	&	0.0	\\\hline
0.3	&	448	&	600	&	3.59	&	315456	&	600	&	0.39	&	1271259	&	600	&	2.05	&	   5	&	5	&	0.0	\\\hline
0.3	&	545	&	600	&	3.90	&	327630	&	600	&	0.68	&	1145088	&	600	&	2.61	&	   3	&	2	&	0.0	\\\hline
0.4	&	494	&	600	&	5.85	&	256185	&	600	&	2.39	&	905790	&	600	&	4.65	&	   3	&	2	&	0.0	\\\hline
0.4	&	541	&	600	&	4.72	&	253367	&	600	&	1.84	&	1030448	&	600	&	3.61	&	  11	&	6	&	0.0	\\\hline
0.4	&	609	&	600	&	7.71	&	153412	&	600	&	5.66	&	806007	&	600	&	6.65	&	  19	&	9	&	0.0	\\\hline
0.4	&	600	&	600	&	7.36	&	198531	&	600	&	4.48	&	839475	&	600	&	6.07	&	 135	&	73	&	0.0	\\\hline
0.4	&	416	&	600	&	6.16	&	314293	&	600	&	0.61	&	971486	&	600	&	5.32	&	 271	&	116	&	0.0	\\\hline
0.4	&	545	&	600	&	7.30	&	247658	&	600	&	3.12	&	1091092	&	600	&	5.23	&	   9	&	4	&	0.0	\\\hline
0.4	&	481	&	600	&	7.07	&	171085	&	600	&	4.71	&	715587	&	600	&	6.06	&	  79	&	34	&	0.0	\\\hline
0.4	&	501	&	600	&	6.05	&	201533	&	600	&	3.00	&	895046	&	600	&	4.78	&	 107	&	51	&	0.0	\\\hline
0.4	&	487	&	600	&	5.05	&	320172	&	600	&	0.43	&	1156723	&	600	&	3.27	&	  61	&	28	&	0.0	\\\hline
0.4	&	436	&	600	&	6.01	&	296422	&	600	&	1.65	&	921533	&	600	&	4.70	&	   5	&	3	&	0.0	\\\hline
0.5	&	527	&	600	&	7.81	&	289631	&	600	&	2.71	&	872174	&	600	&	6.63	&	   9	&	6	&	0.0	\\\hline
0.5	&	540	&	600	&	7.06	&	217761	&	600	&	3.98	&	1005368	&	600	&	5.54	&	   9	&	4	&	0.0	\\\hline
0.5	&	577	&	600	&	10.99	&	181633	&	600	&	8.32	&	756893	&	600	&	9.44	&	  41	&	18	&	0.0	\\\hline
0.5	&	466	&	600	&	10.32	&	185921	&	600	&	6.74	&	828527	&	600	&	8.92	&	 155	&	66	&	0.0	\\\hline
0.5	&	399	&	600	&	8.69	&	285081	&	600	&	4.01	&	1112264	&	600	&	6.33	&	  93	&	42	&	0.0	\\\hline
0.5	&	516	&	600	&	9.95	&	229517	&	600	&	6.27	&	1029917	&	600	&	7.92	&	   5	&	3	&	0.0	\\\hline
0.5	&	532	&	600	&	9.93	&	177439	&	600	&	6.56	&	722884	&	600	&	8.80	&	 167	&	61	&	0.0	\\\hline
0.5	&	492	&	600	&	8.37	&	185644	&	600	&	5.41	&	800682	&	600	&	7.00	&	 137	&	53	&	0.0	\\\hline
0.5	&	413	&	600	&	6.80	&	312503	&	600	&	0.71	&	1072587	&	600	&	4.79	&	 179	&	68	&	0.0	\\\hline
0.5	&	498	&	600	&	9.14	&	265445	&	600	&	3.35	&	921716	&	600	&	7.16	&	  11	&	6	&	0.0	\\\hline
0.6	&	410	&	600	&	9.58	&	267100	&	600	&	4.59	&	943087	&	600	&	8.81	&	  19	&	10	&	0.0	\\\hline
0.6	&	457	&	600	&	9.72	&	235091	&	600	&	5.93	&	982814	&	600	&	7.90	&	   7	&	4	&	0.0	\\\hline
0.6	&	482	&	600	&	13.91	&	185440	&	600	&	11.24	&	803205	&	600	&	12.28	&	  35	&	13	&	0.0	\\\hline
0.6	&	539	&	600	&	13.32	&	210315	&	600	&	10.59	&	953491	&	600	&	11.90	&	 129	&	54	&	0.0	\\\hline
0.6	&	430	&	600	&	11.04	&	284537	&	600	&	5.83	&	1238524	&	600	&	8.07	&	  77	&	27	&	0.0	\\\hline
0.6	&	381	&	600	&	12.34	&	3257574	&	600	&	7.32	&	1161734	&	600	&	10.01	&	  19	&	13	&	0.0	\\\hline
0.6	&	482	&	600	&	11.51	&	210088	&	600	&	8.49	&	782945	&	600	&	11.59	&	  29	&	13	&	0.0	\\\hline
0.6	&	404	&	600	&	11.32	&	4197474	&	600	&	8.19	&	946277	&	600	&	8.59	&	  35	&	7	&	0.0	\\\hline
0.6	&	338	&	600	&	8.61	&	294116	&	600	&	3.60	&	1139723	&	600	&	5.89	&	  83	&	31	&	0.0	\\\hline
0.6	&	392	&	600	&	11.57	&	221902	&	600	&	7.95	&	950386	&	600	&	10.87	&	  25	&	12	&	0.0	\\\hline
0.7	&	388	&	600	&	10.65	&	274783	&	600	&	6.39	&	1111097	&	600	&	10.00	&	  25	&	9	&	0.0	\\\hline
0.7	&	437	&	600	&	12.04	&	1248035	&	600	&	7.93	&	1148012	&	600	&	9.02	&	   9	&	3	&	0.0	\\\hline
0.7	&	465	&	600	&	18.34	&	174297	&	600	&	15.30	&	779013	&	600	&	17.12	&	  31	&	12	&	0.0	\\\hline
0.7	&	512	&	600	&	16.44	&	231018	&	600	&	12.05	&	1175087	&	600	&	13.98	&	 275	&	85	&	0.0	\\\hline
0.7	&	417	&	600	&	13.66	&	273481	&	600	&	7.23	&	1502260	&	600	&	9.33	&	  57	&	14	&	0.0	\\\hline
0.7	&	478	&	600	&	15.70	&	3252290	&	600	&	10.15	&	1265411	&	600	&	13.15	&	  29	&	7	&	0.0	\\\hline
0.7	&	379	&	600	&	13.71	&	248709	&	600	&	8.22	&	901649	&	600	&	14.50	&	  47	&	16	&	0.0	\\\hline
0.7	&	499	&	600	&	14.11	&	4225837	&	600	&	10.31	&	1029597	&	600	&	11.14	&	  21	&	4	&	0.0	\\\hline
0.7	&	435	&	600	&	11.55	&	260515	&	600	&	7.41	&	1248500	&	600	&	7.51	&	  37	&	11	&	0.0	\\\hline
0.7	&	479	&	600	&	15.76	&	218194	&	600	&	12.32	&	1027522	&	600	&	13.59	&	  45	&	19	&	0.0	\\\hline
0.8	&	190	&	600	&	11.31	&	269689	&	600	&	4.56	&	1254085	&	601	&	6.79	&	  13	&	4	&	0.0	\\\hline
0.8	&	349	&	600	&	13.71	&	1249086	&	600	&	9.74	&	1241374	&	600	&	11.32	&	  19	&	5	&	0.0	\\\hline
0.8	&	404	&	600	&	19.46	&	242645	&	600	&	14.33	&	1079730	&	600	&	14.08	&	  43	&	7	&	0.0	\\\hline
0.8	&	350	&	600	&	18.80	&	283263	&	600	&	12.32	&	1400882	&	600	&	14.89	&	  21	&	7	&	0.0	\\\hline
0.8	&	248	&	600	&	12.14	&	301559	&	600	&	4.03	&	1874137	&	600	&	7.84	&	  91	&	15	&	0.0	\\\hline
0.8	&	435	&	600	&	20.38	&	3257129	&	600	&	13.67	&	1360080	&	600	&	14.44	&	  25	&	4	&	0.0	\\\hline
0.8	&	373	&	600	&	15.30	&	231575	&	600	&	11.75	&	1175845	&	600	&	13.51	&	 143	&	37	&	0.0	\\\hline
0.8	&	423	&	600	&	18.67	&	4240437	&	600	&	14.09	&	1286691	&	600	&	11.36	&	  51	&	6	&	0.0	\\\hline
0.8	&	360	&	600	&	14.59	&	285433	&	600	&	7.82	&	1325850	&	600	&	9.56	&	  47	&	10	&	0.0	\\\hline
0.8	&	363	&	600	&	18.74	&	249145	&	600	&	12.84	&	1241172	&	600	&	15.59	&	  93	&	17	&	0.0	\\\hline
0.9	&	129	&	600	&	3.43	&	279839	&	600	&	7.02	&	333869	&	94	&	0.00	&	  27	&	7	&	0.0	\\\hline
0.9	&	247	&	600	&	13.69	&	1294713	&	600	&	3.66	&	2118555	&	600	&	2.69	&	  47	&	6	&	0.0	\\\hline
0.9	&	57	&	239	&	0.10	&	79764	&	170	&	0.0	&	184514	&	58	&	0.00	&	  47	&	6	&	0.0	\\\hline
0.9	&	177	&	600	&	1.60	&	123932	&	209	&	0.0	&	314744	&	80	&	0.00	&	   5	&	2	&	0.0	\\\hline
0.9	&	102	&	334	&	0.10	&	74031	&	115	&	0.0	&	148576	&	38	&	0.00	&	 101	&	8	&	0.0	\\\hline
0.9	&	74	&	297	&	0.10	&	3220333	&	441	&	0.0	&	230267	&	69	&	0.00	&	  55	&	7	&	0.0	\\\hline
0.9	&	112	&	600	&	0.82	&	201689	&	444	&	0.0	&	172017	&	52	&	0.00	&	  63	&	8	&	0.0	\\\hline
0.9	&	295	&	600	&	14.90	&	4282430	&	611	&	3.91	&	1440695	&	424	&	0.00	&	  97	&	11	&	0.0	\\\hline
0.9	&	236	&	600	&	13.04	&	302915	&	600	&	7.21	&	1208680	&	371	&	0.00	&	  85	&	12	&	0.0	\\\hline
0.9	&	131	&	447	&	0.10	&	262283	&	600	&	0.77	&	350219	&	100	&	0.00	&	89	&	10	&	0.0	\\\hline	
\end{tabular}
\end{table}}
{\tiny
\begin{table}[ht]
    \centering
    \begin{tabular}{|l|l|l|l|l|l|l|l|l|l|l|l|l|}\hline
	&	\multicolumn{4}{|c|}{No bound tightening}&	\multicolumn{4}{|c|}{Light bound tightening}&	\multicolumn{4}{|c|}{Strong bound tightening}\\\hline
 $\alpha$ 	& 	nodes 	&          	\#lps 	&        	time 	&  	         GAP\% 	& 	nodes 	&          	\#lps 	&        	time 	&  	         GAP\% & 	nodes 	&          	\#lps 	&        	time 	&  	         GAP\%   	\\\hline
0.3	& 	110411	&          	110414	&        	600	&  	1.38	&  	   45	&  	  2081	&  	  7.73	&  	0.0	&  	   3	&  	 1389	&  	  7.54	&  	0.0	\\\hline
0.3	& 	104983	&          	104986	&        	600	&  	1.37	&  	  117	&  	  5275	&  	  7.21	&  	0.0	&  	   5	&  	 2382	&  	  3.65	&  	0.0	\\\hline
0.3	& 	108667	&          	108670	&        	600	&  	1.15	&  	   75	&  	  3343	&  	  4.83	&  	0.0	&  	   3	&  	 1259	&  	  2.23	&  	0.0	\\\hline
0.3	& 	103993	&          	103996	&        	600	&  	2.66	&  	 1333	&  	 61116	&  	 83.63	&  	0.0	&  	  25	&  	11219	&  	 14.24	&  	0.0	\\\hline
0.3	& 	105677	&          	105680	&        	600	&  	2.29	&  	  911	&  	 41552	&  	 53.47	&  	0.0	&  	  39	&  	13992	&  	 16.67	&  	0.0	\\\hline
0.3	& 	109691	&          	109694	&        	600	&  	2.12	&  	  229	&  	 10158	&  	 13.09	&  	0.0	&  	  15	&  	 4617	&  	  6.15	&  	0.0	\\\hline
0.3	& 	108237	&          	108240	&        	600	&  	2.29	&  	  489	&  	 21752	&  	 29.30	&  	0.0	&  	   9	&  	 4199	&  	  5.65	&  	0.0	\\\hline
0.3	& 	105505	&          	105508	&        	600	&  	2.43	&  	  521	&  	 23854	&  	 31.75	&  	0.0	&  	   5	&  	 3101	&  	  4.37	&  	0.0	\\\hline
0.3	& 	102763	&          	102766	&        	600	&  	1.91	&  	  455	&  	 20895	&  	 27.91	&  	0.0	&  	   5	&  	 3387	&  	  4.67	&  	0.0	\\\hline
0.3	& 	109421	&          	109424	&        	600	&  	1.39	&  	   55	&  	  2501	&  	  3.51	&  	0.0	&  	   3	&  	 1105	&  	  1.84	&  	0.0	\\\hline
0.4	& 	114313	&          	114316	&        	600	&  	1.96	&  	   65	&  	  2953	&  	  3.92	&  	0.0	&  	   3	&  	 1381	&  	  2.49	&  	0.0	\\\hline
0.4	& 	110321	&          	110324	&        	600	&  	2.17	&  	  341	&  	 14440	&  	 19.52	&  	0.0	&  	  11	&  	 3998	&  	  5.79	&  	0.0	\\\hline
0.4	& 	114203	&          	114206	&        	600	&  	1.96	&  	  769	&  	 33128	&  	 44.17	&  	0.0	&  	  19	&  	 6852	&  	  8.79	&  	0.0	\\\hline
0.4	& 	110271	&          	110274	&        	600	&  	4.17	&  	 7603	&  	346771	&  	459.69	&  	0.0	&  	 135	&  	60330	&  	 73.05	&  	0.0	 \\\hline
0.4	& 	110005	&          	110008	&        	600	&  	3.88	&  	10303	&  	463434	&  	600	&  	1.3e-03	&  	 271	&  	99345	&  	115.94	&  	0.0	 \\\hline
0.4	& 	112983	&          	112986	&        	600	&  	2.63	&  	  123	&  	  5224	&  	  7.01	&  	0.0	&  	   9	&  	 2404	&  	  3.58	&  	0.0	\\\hline
0.4	& 	113319	&          	113322	&        	600	&  	3.41	&  	 3113	&  	137531	&  	177.20	&  	0.0	&  	  79	&  	28258	&  	 33.54	&  	0.0	\\\hline
0.4	& 	113269	&          	113272	&        	600	&  	3.77	&  	 4597	&  	206674	&  	270.09	&  	0.0	&  	 107	&  	42825	&  	 51.43	&  	0.0	\\\hline
0.4	& 	108125	&          	108128	&        	600	&  	3.07	&  	 1835	&  	 82792	&  	107.60	&  	0.0	&  	  61	&  	23760	&  	 28.38	&  	0.0	\\\hline
0.4	& 	113003	&          	113006	&        	600	&  	2.28	&  	   85	&  	  3755	&  	  5.10	&  	0.0	&  	   5	&  	 1622	&  	  2.77	&  	0.0	\\\hline
0.5	& 	118889	&          	118892	&        	600	&  	2.72	&  	  233	&  	  9712	&  	 13.19	&  	0.0	&  	   9	&  	 3874	&  	  5.63	&  	0.0	\\\hline
0.5	& 	115619	&          	115622	&        	600	&  	3.02	&  	  207	&  	  8398	&  	 11.43	&  	0.0	&  	   9	&  	 2965	&  	  4.23	&  	0.0	\\\hline
0.5	& 	117685	&          	117688	&        	600	&  	3.18	&  	 1785	&  	 76149	&  	 98.70	&  	0.0	&  	  41	&  	14242	&  	 18.15	&  	0.0	\\\hline
0.5	& 	117515	&          	117518	&        	600	&  	5.66	&  	 3939	&  	174038	&  	220.90	&  	0.0	&  	 155	&  	55206	&  	 65.84	&  	0.0	\\\hline
0.5	& 	116919	&          	116922	&        	600	&  	5.70	&  	 3231	&  	141810	&  	177.74	&  	0.0	&  	  93	&  	35390	&  	 41.94	&  	0.0	\\\hline
0.5	& 	120505	&          	120508	&        	600	&  	3.34	&  	  123	&  	  4958	&  	  6.74	&  	0.0	&  	   5	&  	 1925	&  	  3.04	&  	0.0	\\\hline
0.5	& 	117217	&          	117220	&        	600	&  	4.81	&  	 3369	&  	139806	&  	185.61	&  	0.0	&  	 167	&  	48601	&  	 60.57	&  	0.0	\\\hline
0.5	& 	117795	&          	117798	&        	600	&  	5.42	&  	 2623	&  	110557	&  	144.71	&  	0.0	&  	 137	&  	42937	&  	 52.99	&  	0.0	\\\hline
0.5	& 	114831	&          	114834	&        	600	&  	4.30	&  	 3717	&  	162722	&  	210.93	&  	0.0	&  	 179	&  	57501	&  	 68.47	&  	0.0	\\\hline
0.5	& 	118681	&          	118684	&        	600	&  	3.75	&  	  271	&  	 11069	&  	 14.94	&  	0.0	&  	  11	&  	 4045	&  	  5.75	&  	0.0	\\\hline
0.6	& 	122649	&          	122652	&        	600	&  	3.79	&  	  343	&  	 13535	&  	 18.47	&  	0.0	&  	  19	&  	 6893	&  	  9.60	&  	0.0	\\\hline
0.6	& 	120833	&          	120836	&        	600	&  	3.83	&  	  113	&  	  4353	&  	  6.06	&  	0.0	&  	   7	&  	 2361	&  	  3.88	&  	0.0	\\\hline
0.6	& 	122259	&          	122262	&        	600	&  	4.10	&  	  961	&  	 40108	&  	 50.37	&  	0.0	&  	  35	&  	10673	&  	 13.14	&  	0.0	\\\hline
0.6	& 	112075	&          	112078	&        	600	&  	8.07	&  	 4449	&  	184494	&  	235.41	&  	0.0	&  	 129	&  	45165	&  	 54.31	&  	0.0	\\\hline
0.6	& 	114885	&          	114888	&        	600	&  	7.75	&  	 1931	&  	 81916	&  	 98.60	&  	0.0	&  	  77	&  	23555	&  	 27.16	&  	0.0	\\\hline
0.6	& 	117107	&          	117110	&        	600	&  	4.89	&  	  505	&  	 18465	&  	 25.44	&  	0.0	&  	  19	&  	 9271	&  	 13.02	&  	0.0	\\\hline
0.6	& 	117699	&          	117702	&        	600	&  	0.05	&  	  885	&  	 35256	&  	 44.75	&  	0.0	&  	  29	&  	10215	&  	 13.19	&  	0.0	\\\hline
0.6	& 	122133	&          	122136	&        	600	&  	5.69	&  	  227	&  	  9312	&  	 11.05	&  	0.0	&  	  35	&  	 6269	&  	  7.24	&  	0.0	\\\hline
0.6	& 	116661	&          	116664	&        	600	&  	5.50	&  	 1081	&  	 46073	&  	 57.46	&  	0.0	&  	  83	&  	25607	&  	 31.10	&  	0.0	\\\hline
0.6	& 	113123	&          	113126	&        	600	&  	5.65	&  	  573	&  	 22301	&  	 29.69	&  	0.0	&  	  25	&  	 8816	&  	 11.96	&  	0.0	\\\hline
0.7	& 	118717	&          	118720	&        	600	&  	5.33	&  	  405	&  	 15074	&  	 20.44	&  	0.0	&  	  25	&  	 6485	&  	  8.75	&  	0.0	\\\hline
0.7	& 	125675	&          	125678	&        	600	&  	3.40	&  	   65	&  	  2426	&  	  3.32	&  	0.0	&  	   9	&  	 1992	&  	  2.81	&  	0.0	\\\hline
0.7	& 	127137	&          	127140	&        	600	&  	4.47	&  	  341	&  	 13609	&  	 16.54	&  	0.0	&  	  31	&  	 9504	&  	 11.91	&  	0.0	\\\hline
0.7	& 	125293	&          	125296	&        	600	&  	10.91	&  	 3671	&  	140527	&  	179.20	&  	0.0	&  	 275	&  	70938	&  	 85.43	&  	0.0	\\\hline
0.7	& 	126751	&          	126754	&        	600	&  	9.02	&  	  821	&  	 31203	&  	 37.62	&  	0.0	&  	  57	&  	13047	&  	 14.34	&  	0.0	\\\hline
0.7	& 	128313	&          	128316	&        	600	&  	4.95	&  	  255	&  	  8927	&  	 11.85	&  	0.0	&  	  29	&  	 5343	&  	  6.72	&  	0.0	\\\hline
0.7	& 	128455	&          	128458	&        	600.02	&  	6.75	&  	  805	&  	 29345	&  	 38.77	&  	0.0	&  	  47	&  	12328	&  	 15.88	&  	0.0	\\\hline
0.7	& 	86603	&          	86606	&        	600	&  	6.19	&  	  133	&  	  5124	&  	  6.12	&  	0.0	&  	  21	&  	 3327	&  	  4.02	&  	0.0	\\\hline
0.7	& 	85403	&          	85406	&        	600	&  	6.42	&  	  307	&  	 12083	&  	 14.88	&  	0.0	&  	  37	&  	 9075	&  	 10.66	&  	0.0	\\\hline
0.7	& 	86247	&          	86250	&        	600	&  	8.57	&  	 1335	&  	 47897	&  	 65.73	&  	0.0	&  	  45	&  	14010	&  	 19.21	&  	0.0	\\\hline
0.8	& 	87179	&          	87182	&        	600	&  	6.49	&  	  147	&  	  5264	&  	  6.91	&  	0.0	&  	  13	&  	 2828	&  	  3.63	&  	0.0	\\\hline
0.8	& 	92415	&          	92418	&        	600	&  	2.64	&  	   63	&  	  2326	&  	  2.94	&  	0.0	&  	  19	&  	 3719	&  	  4.93	&  	0.0	\\\hline
0.8	& 	92645	&          	92648	&        	600	&  	3.22	&  	  187	&  	  6850	&  	  8.20	&  	0.0	&  	  43	&  	 5984	&  	  6.88	&  	0.0	\\\hline
0.8	& 	90021	&          	90024	&        	600.02	&  	11.63	&  	  233	&  	  8261	&  	 10.29	&  	0.0	&  	  21	&  	 5408	&  	  6.98	&  	0.0	\\\hline
0.8	& 	89669	&          	89672	&        	600	&  	9.09	&  	  515	&  	 17362	&  	 21.23	&  	0.0	&  	  91	&  	13655	&  	 15.06	&  	0.0	\\\hline
0.8	& 	91903	&          	91906	&        	600	&  	4.36	&  	  131	&  	  4484	&  	  5.79	&  	0.0	&  	  25	&  	 3540	&  	  4.34	&  	0.0	\\\hline
0.8	& 	91247	&          	91250	&        	600	&  	9.23	&  	  827	&  	 26964	&  	 35.37	&  	0.0	&  	 143	&  	28591	&  	 36.72	&  	0.0	\\\hline
0.8	& 	90817	&          	90820	&        	600	&  	5.39	&  	  141	&  	  5029	&  	  5.95	&  	0.0	&  	  51	&  	 5923	&  	  6.49	&  	0.0	\\\hline
0.8	& 	89261	&          	89264	&        	600	&  	6.20	&  	  321	&  	 11617	&  	 14.22	&  	0.0	&  	  47	&  	 8291	&  	  9.83	&  	0.0	\\\hline
0.8	& 	91311	&          	91314	&        	600	&  	9.41	&  	  505	&  	 17849	&  	 21.37	&  	0.0	&  	  93	&  	14531	&  	 17.13	&  	0.0	\\\hline
0.9	& 	88037	&          	88040	&        	600	&  	6.67	&  	  165	&  	  5300	&  	  6.88	&  	0.0	&  	  27	&  	 4879	&  	  6.54	&  	0.0	\\\hline
0.9	& 	105609	&          	105612	&        	600	&  	0.16	&  	  189	&  	  5922	&  	  7.50	&  	0.0	&  	  47	&  	 4956	&  	  5.58	&  	0.0	\\\hline
0.9	& 	72451	&          	72454	&        	399.49	&  	0.00	&  	  265	&  	  8503	&  	 10.78	&  	0.0	&  	  47	&  	 5669	&  	  6.33	&  	0.0	\\\hline
0.9	& 	93485	&          	93488	&        	600	&  	5.33	&  	   35	&  	  1229	&  	  1.74	&  	0.0	&  	   5	&  	  999	&  	  1.56	&  	0.0	\\\hline
0.9	& 	93089	&          	93092	&        	553.74	&  	0.00	&  	  177	&  	  5565	&  	  7.06	&  	0.0	&  	 101	&  	 7728	&  	  8.37	&  	0.0	\\\hline
0.9	& 	97929	&          	97932	&        	600	&  	2.60	&  	  175	&  	  5352	&  	  7.02	&  	0.0	&  	  55	&  	 6067	&  	  7.05	&  	0.0	\\\hline
0.9	& 	74985	&          	74988	&        	439.06	&  	0.00	&  	  185	&  	  5849	&  	  7.42	&  	0.0	&  	  63	&  	 7175	&  	  8.22	&  	0.0	\\\hline
0.9	& 	81611	&          	81614	&        	473.46	&  	0.00	&  	  165	&  	  5212	&  	  6.57	&  	0.0	&  	  97	&  	 9225	&  	 10.84	&  	0.0	\\\hline
0.9	& 	86963	&          	86966	&        	474.28	&  	0.00	&  	  265	&  	  8471	&  	 10.81	&  	0.0	&  	  85	&  	 9562	&  	 11.54	&  	0.0	\\\hline
0.9	& 	100867	&          	100870	&        	571.30	&  	0.00	&  	  211	&  	  6699	&  	  8.55	&  	0.0	&  	  89	&  	 8512	&  	  9.59	&  	0.0	\\\hline
\end{tabular}
\end{table}}



\end{document}